\documentclass[titlepage, sn-mathphys]{sn-jnl}

\usepackage{amsmath}
\usepackage[labelfont=bf]{subcaption}
\usepackage[labelsep = space]{caption}

\jyear{2022}%

\theoremstyle{thmstyleone}%
\theoremstyle{thmstyletwo}%

\theoremstyle{thmstylethree}%
\DeclareMathOperator*{\argmin}{arg\,min}

\newcommand{\Pop}{\mathit{Pop}}
\newcommand{\Dev}{\mathit{Dev}}
\newcommand{\LI}{\mathit{LI}}
\newcommand{\MI}{\mathit{MI}}

\newcommand{\FLIP}{\mathit{FLIP}}
\newcommand{\SWAP}{\mathit{SWAP}}
\newcommand{\CMB}{\mathit{CMB}}
\newcommand{\NH}{\mathit{NH}}
\newcommand{\NF}{\mathit{NF}}

\begin{document}

\begin{titlepage}

\title[A $k$-medoids Approach to Exploring Districting Plans]{A $k$-medoids Approach to Exploring Districting Plans}

\author*[1]{\fnm{Jared} \sur{Grove}}\email{jared-grove@uiowa.edu}

\author[2]{\pfx{Dr.} \fnm{Suely P.} \sur{Oliveira}}\email{suely-oliveira@uiowa.edu}

\author[3]{\fnm{Anthony} \sur{Pizzimenti}}\email{anthony@mggg.org}

\author[1]{\pfx{Dr.} \fnm{David} \sur{Stewart}}\email{david-e-stewart@uiowa.edu}

\affil[1]{\orgdiv{Mathematics}, \orgname{University of Iowa}, \orgaddress{\street{2 West Washington St}, \city{Iowa City}, \postcode{52242}, \state{Iowa}, \country{United States}}}

\affil[2]{\orgdiv{Computer Science}, \orgname{University of Iowa}, \orgaddress{\street{2 West Washington St}, \city{Iowa City}, \postcode{52242}, \state{IA}, \country{United States}}}

\affil[3]{\orgdiv{Metric Geometry and Gerrymandering Group}, \orgname{Tufts University}, \orgaddress{\street{10 Upper Campus Rd}, \city{Medford}, \postcode{02155}, \state{MA}, \country{United States}}}

\abstract{ Researchers and legislators alike continue the search for methods of drawing fair districting plans. A districting plan is a partition of a state's subdivisions (e.g. counties, voting precincts, etc.). By modeling these districting plans as graphs, they are easier to create, store, and operate on. Since graph partitioning with balancing populations is a computationally intractable (NP-hard) problem most attempts to solve them use heuristic methods. In this paper, we present a variant on the $k$-medoids algorithm where, given a set of initial medoids, we find a partition of the graph's vertices to admit a districting plan. We first use the $k$-medoids process to find an initial districting plan, then use local search strategies to fine-tune the results, such as reducing population imbalances between districts. We also experiment with coarsening the graph to work with fewer vertices. The algorithm is tested on Iowa and Florida using 2010 census data to evaluate the effectiveness.}


\maketitle

\vspace{-1.17cm}
\begin{center}
\textbf{Statements and Declarations}\\
We wish to acknowledge the support of the National Science Foundation for Anthony Pizzimenti through their grant NSF-1407216.
\end{center}

\end{titlepage}

\section{Introduction}
\label{section:Intro}
    \subsection{Districting Plans}
    \label{section:IntroDP}
    
    Every $10$ years the United States goes through a reapportionment process, where states may gain or lose seats in the House of Representatives based on population changes. States with more than one representative may need to redistrict, or redraw the lines for their congressional districts based on how the populations have shifted within the state or because the number of representatives for that state changed.
    
\begin{figure}
	\subcaptionbox{\label{fig:IowaMap}}{\centering \includegraphics[width = 0.48\textwidth]{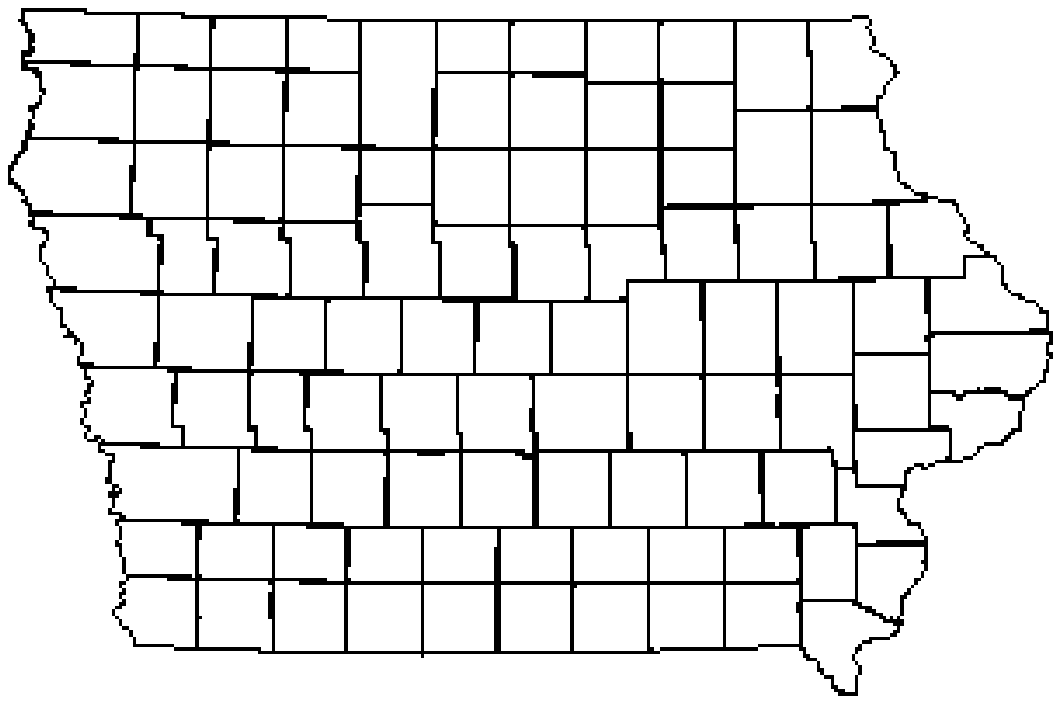}} 
	\subcaptionbox{\label{fig:IowaGraph}}{\centering \includegraphics[width = 0.48\textwidth]{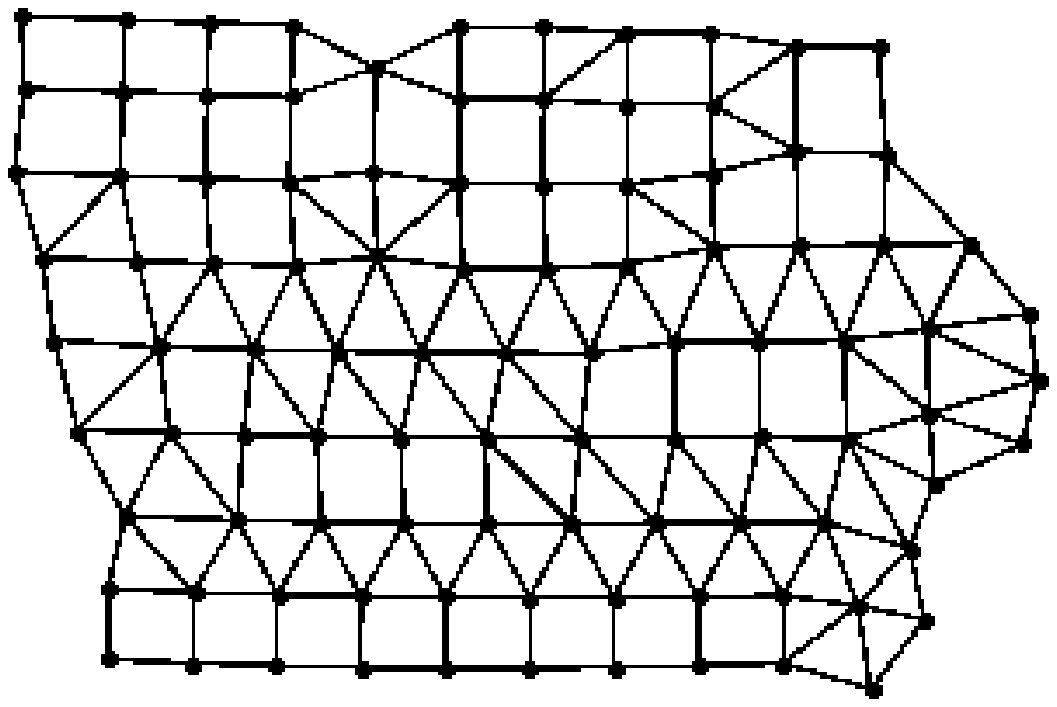}} 
	\caption{County map of Iowa given in (\ref{fig:IowaMap}) and graphical representation of Iowa given in (\ref{fig:IowaGraph})}
	\label{fig:IowaMaps}
\end{figure}

    When a state redistricts, lines must be drawn to divide the state into $k$ districts, one for each Congressional seat. This can be represented mathematically by using a graph or network structure, where the state is a graph $G = (V,E)$ where $V$ is the set of vertices and $E$ is the set of edges. Each state is composed of many voter tabulation districts (VTDs) that need to be assigned to districts, and some states require that larger existing political structures remain intact. For example, Iowa requires that counties are not split apart. We will refer to the basic building blocks of a districting plan as tabulation blocks (TBs), which can represent counties, VTDs, precincts, or any combination of them. If there are $N$ TBs, they will be represented in $G$ as vertices $\left\{ \, v_i \mid i = 1, \ldots, N \,\right\}= V$. The edges in $G$, $\left \{ v_{i},\, v_{j} \right\} \in E$, will represent a shared geographical border between vertices $v_i$ and $v_j$. As an example of county level TBs, the state map and graphical representation of Iowa can be seen in Figure~\ref{fig:IowaMaps}.

    To build a districting plan, each vertex $v_j$ needs to be assigned to a district $v_j \in D_i$ in the districting plan $\left\{ \, D_i \mid i = 1,\ldots, k \, \right\}$ which is a partition of the TBs. 
    The TBs must be assigned to districts such that the districts fit defined criteria for contiguity, compactness, and population equality.
    The criteria for contiguity requires a district to be a single connected region without any holes. The criteria for compactness is less strictly defined, in that a district is considered compact when it covers the smallest possible geographic region. This is often interpreted to mean that a district should be as close to square or circular as possible. 
    A variety of compactness measures have been used previously, including Polsby-Popper \cite{Polsby}, Reock \cite{Reock}, Schwartzberg \cite{Schwartz}, Convex Hull \cite{RestOfComp}, X-Symmetry \cite{RestOfComp}, and Length-Width Ratio \cite{LWR}. All listed measures return a compactness value between $0$ and $1$ with $1$ being considered compact and $0$ being not compact. Non-compact districts are often long and winding with a bizarre shape \cite{Bizarre}, and are often rejected because non-compactness is a hallmark of gerrymandered districts. In this work we will be using the Convex Hull as our compactness measure as it has been used in many court cases \cite{chen_rodden_2015} and, as shown in \cite{Article4-Jumble}, is less affected by choice of map projection. The Convex Hull measure is defined by the ratio of district area to the area of the minimum convex polygon bounding the district. We will report the compactness value as the mean of the compactness values of each district in the districting plan. Lastly, population equality requires that each district has approximately equal populations. 
    Thus an ideal district would have a population equal to the total population of the state divided by the number of districts $k$:
    \begin{equation}
        \Pop^* = \frac{1}{k} \sum_{i=1}^{N} \Pop(v_i), 
    \end{equation}
    where $\Pop(v_i)$ is the population of the TB represented by vertex $v_i$. 
    In order to evaluate a districting plan we will compare the population of each district, $\Pop(D_i) = \sum_{j : v_j \in D_i} \Pop(v_j)$ to $\Pop^*$. The deviation of a districting plan will be the maximum absolute percent deviation of the districts:
    \begin{equation}
	\Dev = \max_{i \in \{1, \ldots, k\}} \left( \frac{\lvert \Pop(D_i) - \Pop^* \rvert}{\Pop^*} \cdot 100\% \right)
   \end{equation}
    We aim to create an algorithm to generate districting plans with  $\Dev<1$, which is equivalent to a districting plan that has maximum population deviation of less than $1\%$. We will do this only using the geographic data and population information.

    \subsection{Gerrymandering}
    \label{section:IntroGerry}
    
    While redistricting is a necessary process that can allow the people to be properly represented in Congress, it can also be manipulated by the agents that draw the lines. 
    This is known as gerrymandering and can be done to retain political power, take power away from opponents, or to further various other political goals. 
    In the worst cases, political officials are essentially able to select their voters, rather than voter selecting their officials. 
    
    The most obvious gerrymanders are often identified by the bizarreness of their shape \cite{Bizarre}, for example long narrow districts that wind throughout the state. 
    These types of districts are often easily identifiable using compactness measures. 
    However, some gerrymanders are much less obvious to the eye and can be much harder to identify. 
    To avoid this problem various algorithms have been proposed to create districting plans, effectively removing politicians and their agents from the redistricting process. However, algorithm based districting does not necessarily create politically neutral districts as the geographic and population data can correlate with political affiliation.
    
    Generating districting plans using methods that reduce political and other bias can also help identify gerrymandering  \cite{EnsambleOutliers}.
    A gerrymandered districting plan would have political outcomes (that is, election outcomes) well outside the range of election outcomes for districting plans created without political bias.
    Generating large numbers of districting plans without political bias, and then determining the expected election outcomes for these generated  plans shows the range of outcomes to be expected where there is no gerrymander.
    Expected election outcomes for a given districting plan well outside this range then becomes strong evidence of a gerrymander
    
    \subsection{Other Approaches}
    \label{section:IntoOA}
    
    The majority of computer-based districting approaches are either sampling or optimization based methods.

    This is because many of the formulations for this type of problem have been shown to be NP-Complete (Subset Sum \cite{computersAndIntractibility}).
    
    Sampling based methods are generally concerned with generating many plans quickly. Often times these are flood fill or multi-kernel methods that start by randomly selecting several initial centers to iteratively build districts from. Typically these will completely build one district, then the second, and so on until there are the proper number of districts. This process has been discussed and implemented many times \cite{vickrey_1961-Sampling,GROUP-Sampling,SouthCarolina-Sampling,thoreson_liittschwager_1967-Sampling}. A potential issue with this method is that it can create enclaves, groups of TBs that cannot be added to any district, and implementations have to monitor this and will restart if enclaves appear. A further variation of the basic sampling method is mentioned in \cite{Article11-Algorithms}, growing the districts in simultaneously instead of sequentially. This should keep the shape of the districts more consistent compared to the sequential version. A similar method \cite{chen_rodden_2015}, considers all TBs as their own district and combines them until the correct number of districts has been achieved. 
    
    The second family of approaches is based on optimization, and seeks to find the best plan under some objective function. These can be roughly subdivided into two more categories; clustering and heuristics. 
    
    Clustering algorithms typically follow a $k$-means or warehouse location approach to assign TBs to cluster centers while minimizing the sum of squares distances from the TBs to the assigned cluster center. Several approaches using these techniques can be found, \cite{weaver_hess_1963-KMEANS,guest_kanayet_love_2019-KMEANS,cohen-addad_klein_young_2018-KMEANS,Article1-Heuristic}. 
    When following a clustering approach special care needs to be taken to ensure that the dividing lines do not break up the TBs that make up the state; some of these legal technicalities are ignored by certain implementations of these methods (see, for example, \cite{guest_kanayet_love_2019-KMEANS,cohen-addad_klein_young_2018-KMEANS}).
    
    Heuristic algorithms often involve using a mixed integer linear program to model districting and find an optimal solution according to some criteria, often minimizing population deviation. 
    They can be solved with various local searches \cite{nagel_1965-LS,kaiser_1966-LS,ricca_simeone_2008-LS}, tabu search \cite{BOZKAYA200312-TABU,ricca_simeone_2008-LS}, genetic algorithms \cite{PEAR}, simulated annealing \cite{ricca_simeone_2008-LS,fifield_higgins_imai_tarr_2020-SA}, and many more heuristic approaches. 
    While these can do a great job of finding plans according to their objective function, there is always the issue of agreeing on what the plans should optimize for: should they be focused on minimizing cut edges in the graph, getting the populations exactly equal, or some other criteria? 
    
    A third type of algorithm that seems to have less representation in the literature uses random walks over the set of districting plans \cite{fifield_higgins_imai_tarr_2020-SA,Gerrychain}. 
    Random walks move from one districting plan to another in an attempt to understand what the solution space of districting plans looks like.
    These types of methods need many good quality initial seeds in order to effective sample the solution space. 
    Combining them with either the sampling or optimization methods can provide these random walk methods with the starting points they need to explore the space \cite{Gerrychain}.
    
    \subsection{Our Approach}
    \label{section:IntroUS}
    We use a $k$-medoids algorithm to approach this problem. 
    It is a sampling method that borrows ideas from some of the optimization methods. 
    We begin from the graph representation of the state and follow the multikernel procedure mentioned in Section~\ref{section:IntoOA}.
    However, instead of following the sequential building pattern used by most implementation we loosely follow an alternating building pattern mentioned in \cite{Article11-Algorithms}. 
    The basics of the $k$-medoids approach is to start with $k$ TBs as initial district centers, then alternate between assigning TBs to the nearest centers and computing the centers of districts. 
    It is similar to the $k$-means method for clustering \cite{weaver_hess_1963-KMEANS,guest_kanayet_love_2019-KMEANS,cohen-addad_klein_young_2018-KMEANS,Article1-Heuristic}.
    
    However, instead of using the physical distances between units, we use graph-based distances. 
    Furthermore, in $k$-medoids, the centers, or medoids, are required to be vertices on the graph as opposed to any point in space. 
    Following the $k$-medoids process we move into a local search phase. 
    We use the traditional local search strategies, flip and swap neighborhoods, along with a tabu criteria to fine-tune solutions. 
    In states with many TBs we also implement a coarsening strategy to work with fewer TBs, then follow an uncoarsening schedule to allow for more fine-tuning as better solutions are found. 
    To the authors' knowledge, the first publication of this coarsening/uncoarsening strategy for districting plans is Magkeby and Mosessoon \cite{magleby_mosesson_2018-Coarsening}. 
    We explore the implemented algorithms in more detail in Section~\ref{section:Algorithms}.

\section{Algorithms}
\label{section:Algorithms}
    The base of the algorithm described here is a $k$-medoids algorithm, which alternates between assigning points to clusters based on proximity to a medoid and computing the medoids of the clusters. However, in the most basic form, $k$-medoids algorithms fail to create districts with equal populations. Several additional steps are added to the basic $k$-medoids framework to account for this issue and improve overall performance. Before describing the algorithm as a whole, it is first useful to walk through these additional steps. These include various forms of local search as well as graph coarsening and uncoarsening.
    
    \subsection{Local Search}
    With the $k$-medoids process used here, deviation does not strictly decrease, so finding a districting plan with a low deviation is not guaranteed.
    
    To increase the likelihood of finding a good solution we look to fine tune the output of the $k$-medoids algorithm by making small changes to the resultant districting plans. These small changes can be thought of as finding neighboring districting plans. Methods of identifying neighboring districting plans are known as neighborhood functions.
    
    We consider three types of neighborhood functions -  Flip, Swap, and a Combination Search (CMB). The Flip neighborhood function generates all districting plans that have exactly one TB assigned to a different district. The Swap neighborhood function generates all districting plans that have exactly two TBs that have mutually changed districts. The CMB neighborhood function generates all districting plans that are in either the Flip neighborhood or Swap neighborhood. In general, we will abbreviate an arbitrary neighborhood function as $\NF$.
    
    An issue that arises with these neighborhood functions is that the number of neighboring districting plans generated can be very large, making it infeasible to enumerate them and find the neighboring districting plan that would have the lowest deviation. 
    
    To get around this, we only generate neighboring districting plans according to specific criteria. From the initial districting plan we consider all pairs of neighboring districts $\left\{ D_m,D_\ell \right\}$ with $\Pop(D_m) > \Pop(D_\ell)$. From these we find the pair of districts with the largest population disparity. We use these districts differently for each method:

\begin{itemize}
        \item[Flip:] Find the set of all vertices $v_i\in D_m$ that share a geographical border with any vertex in $D_\ell$ and that are not articulation points of the graph of $D_m$; articulation points are vertices that when removed from a graph increase the number of connected components. This will be stored as \newline $\FLIP = \left\{\,(v_i, v_{-1})\,\right\}$ where $v_{-1}$ is a dummy vertex with $\Pop(v_{-1}) = 0$.
        
        \item[Swap:] Find the set of all pairs of vertices $\SWAP = \left\{\,(v_i, v_j) \mid v_i \in D_m,\, v_j \in D_\ell \, \right\}$ such that neither $v_i$ nor $v_j$ are articulation points of their respective districts, $v_i$ shares a geographical border with a vertex in $D_\ell \backslash\{v_j \}$, and $v_j$ shares a geographical border with a vertex in $D_m \backslash \{v_i\}$.
        
        \item[CMB:] $\CMB = \FLIP \, \cup \, \SWAP$. The set of pairs of vertices that satisfy either Flip or Swap.
    \end{itemize}
    
    The sets of pairs of vertices derived from the above rules represent the districting plans generated according to the corresponding neighborhood function, which we will denote $\NH$. If no districting plans were generated, we select the pair of neighboring districts with the next highest population disparity until there is at least one districting plan generated.
    
    Once a districting plan is generated, the local search considers all generated districting plans and selects the one with the minimum population disparity. This is done by finding

\begin{align}
\begin{split}
\label{LocalSearchMinimizationEquation}
 (v^*_i, v^*_j) = \argmin_{(v_i, v_j)\in \NH} \max  \left( \, \left\lvert \, \frac{[\Pop(D_m) - \Pop(v_j) + \Pop(v_j)] - \Pop^*}{\Pop^*} \, \right\rvert , \right.\\
 \left. \left\lvert \, \frac{[\Pop(D_\ell) + \Pop(v_i) - \Pop(v_j)] - \Pop^*}{\Pop^*}\, \right\rvert \, \right)
 \end{split}
\end{align}
    
    This pair of vertices represents the districting plan with the minimum population disparity. 
    
\begin{figure}
	\subcaptionbox{\label{fig:Flip1}}{\centering \includegraphics[width = 0.48\textwidth]{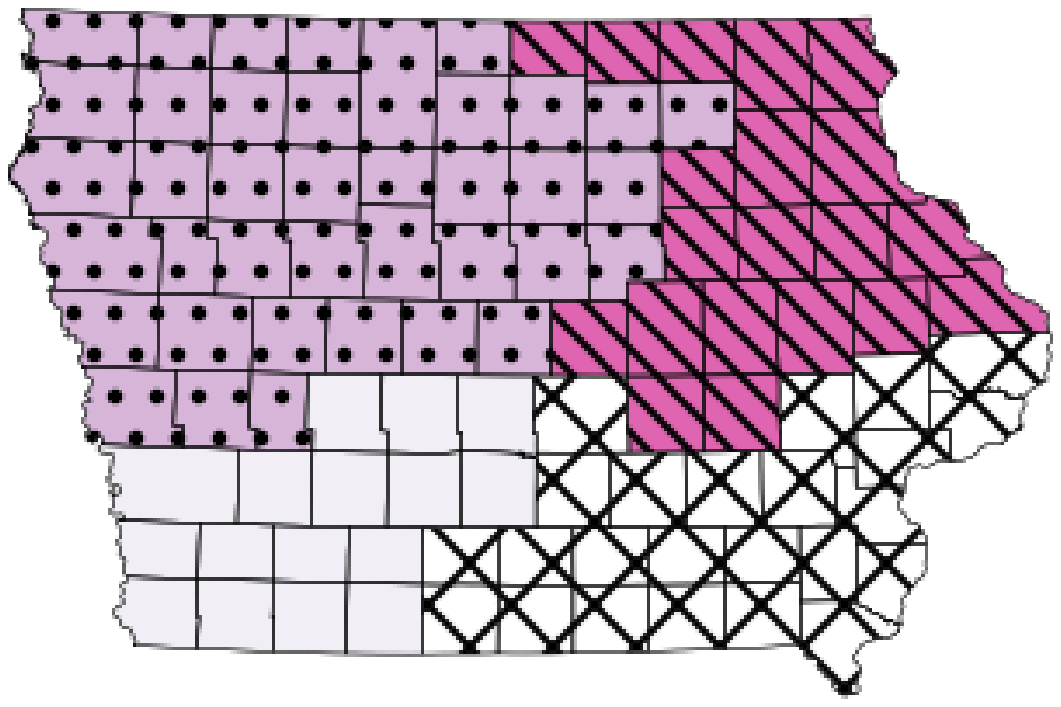}} 
	\subcaptionbox{\label{fig:Flip2}}{\centering \includegraphics[width = 0.48\textwidth]{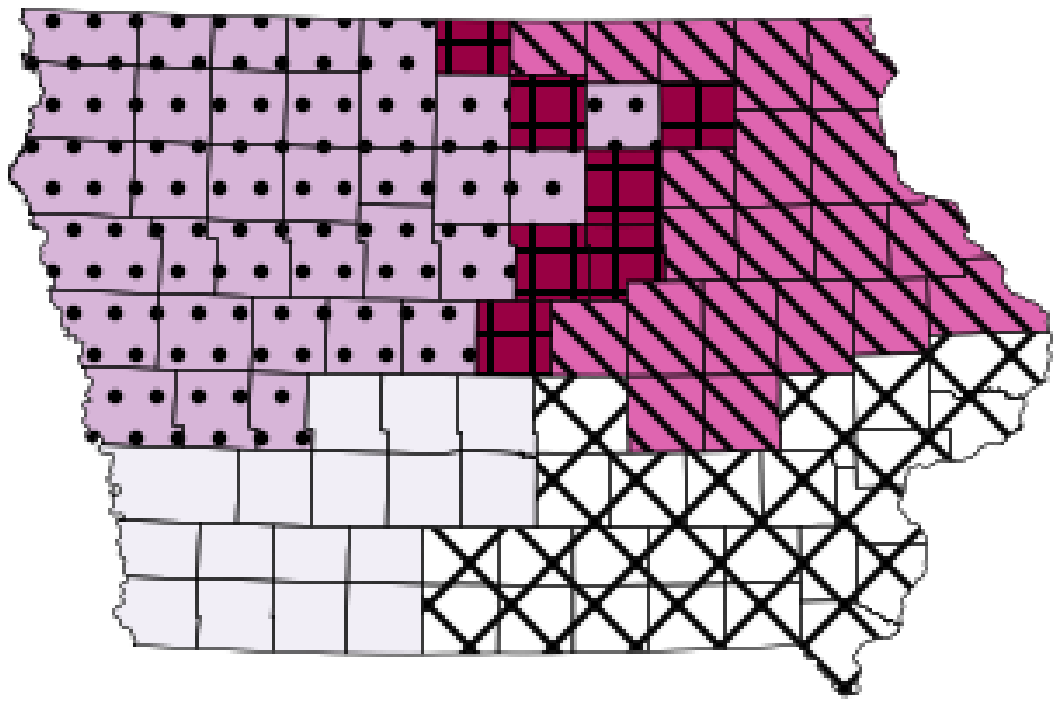}} 
	\subcaptionbox{\label{fig:Flip3}}{\centering \includegraphics[width = 0.48\textwidth]{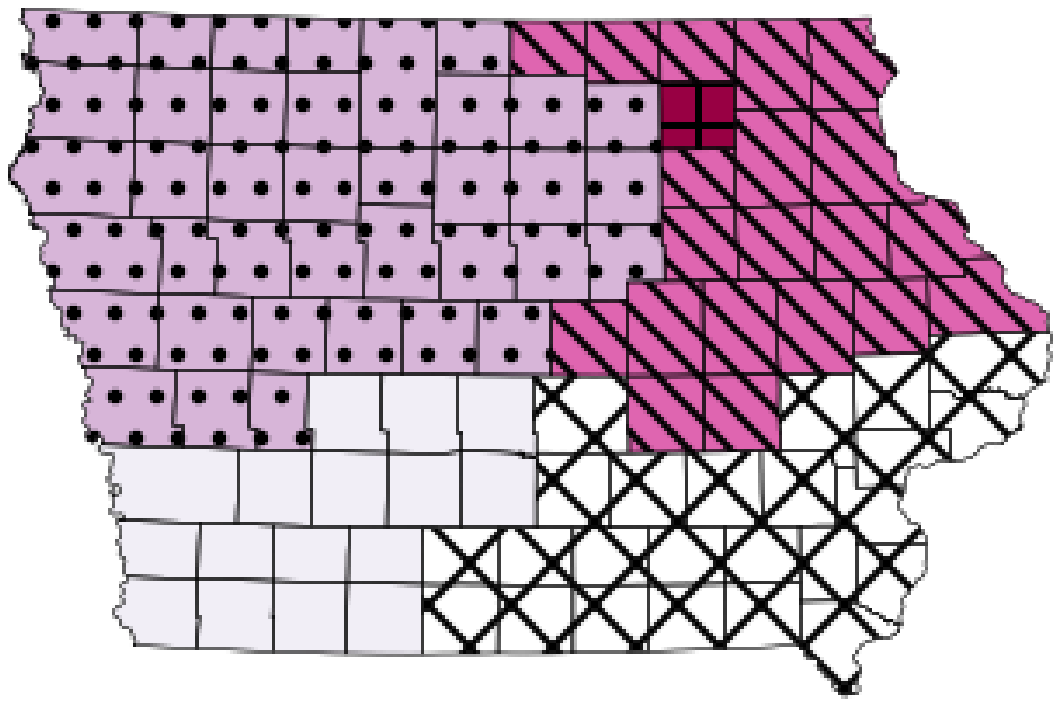}} 
	\subcaptionbox{\label{fig:Flip4}}{\centering \includegraphics[width = 0.48\textwidth]{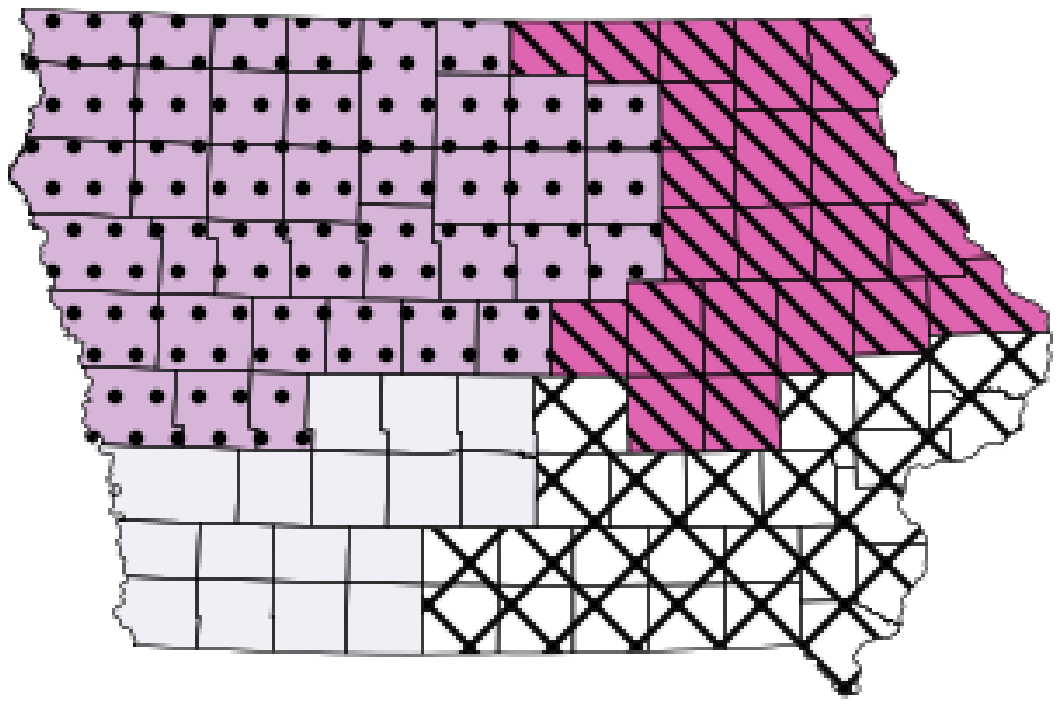}} 
	\caption{Example of Flip neighborhood function. First select Districts that Flip will be performed on, $D_m$: \protect\includegraphics[height = 0.25cm]{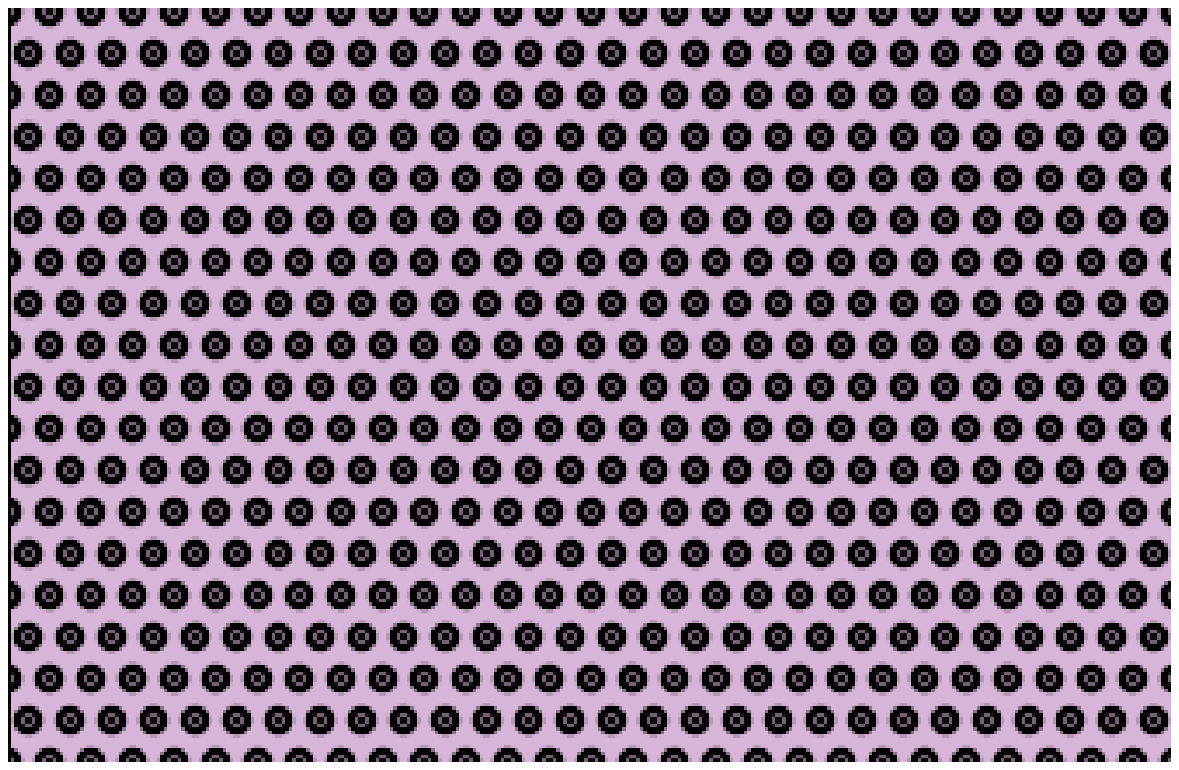} and $D_\ell$: \protect\includegraphics[height =0.25cm]{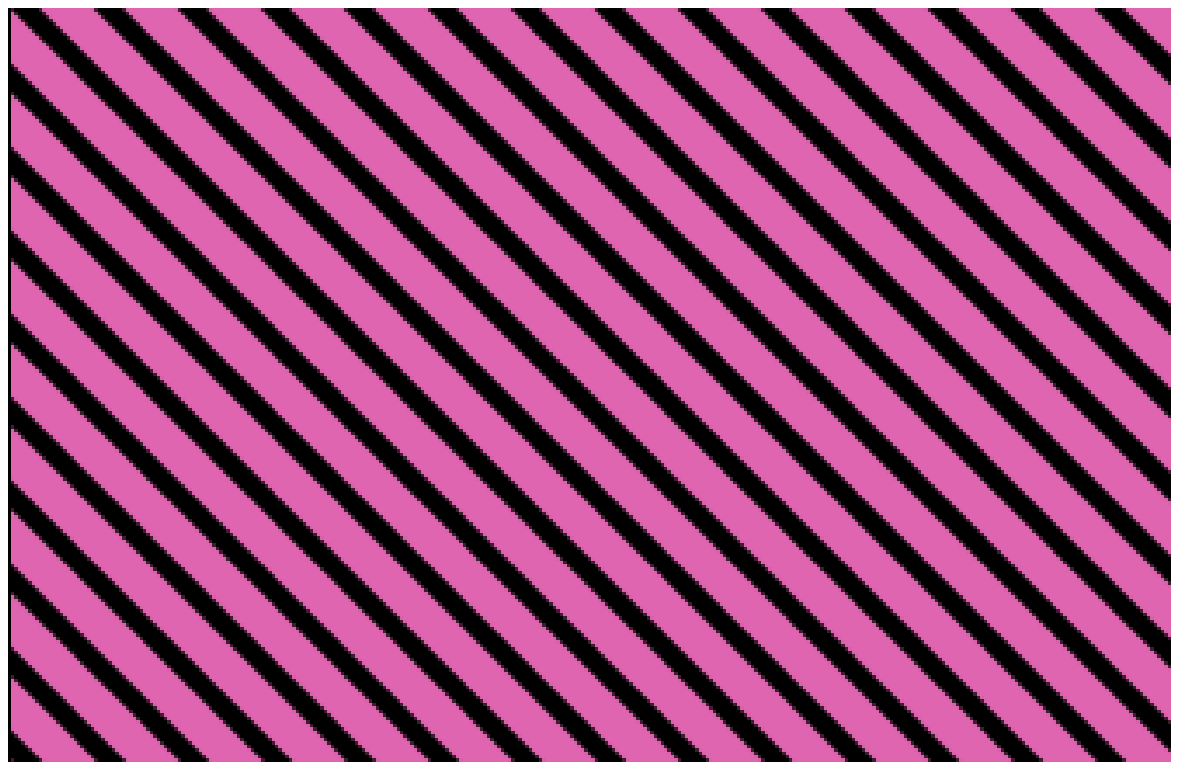} (\ref{fig:Flip1}), then find candidate vertices in $D_m$: \protect\includegraphics[height = 0.25cm]{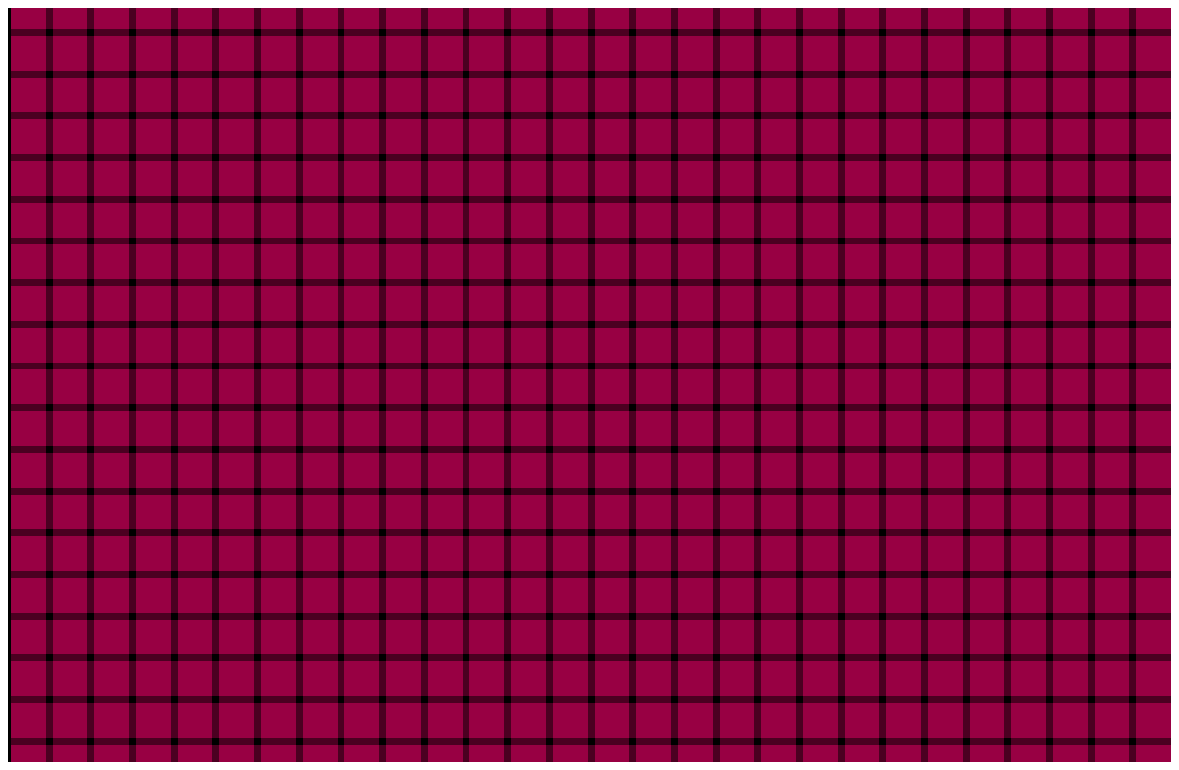} that can be flipped (\ref{fig:Flip2}). Next find candidate vertex that minimizes Equation (\ref{LocalSearchMinimizationEquation})  (\ref{fig:Flip3}) and Flip minimizing vertex from $D_m$ to $D_\ell$ (\ref{fig:Flip4})
	}
	\label{fig:FlipExample}
\end{figure}

\begin{figure}
	\subcaptionbox{\label{fig:Swap1}}{\centering \includegraphics[width = 0.48\textwidth]{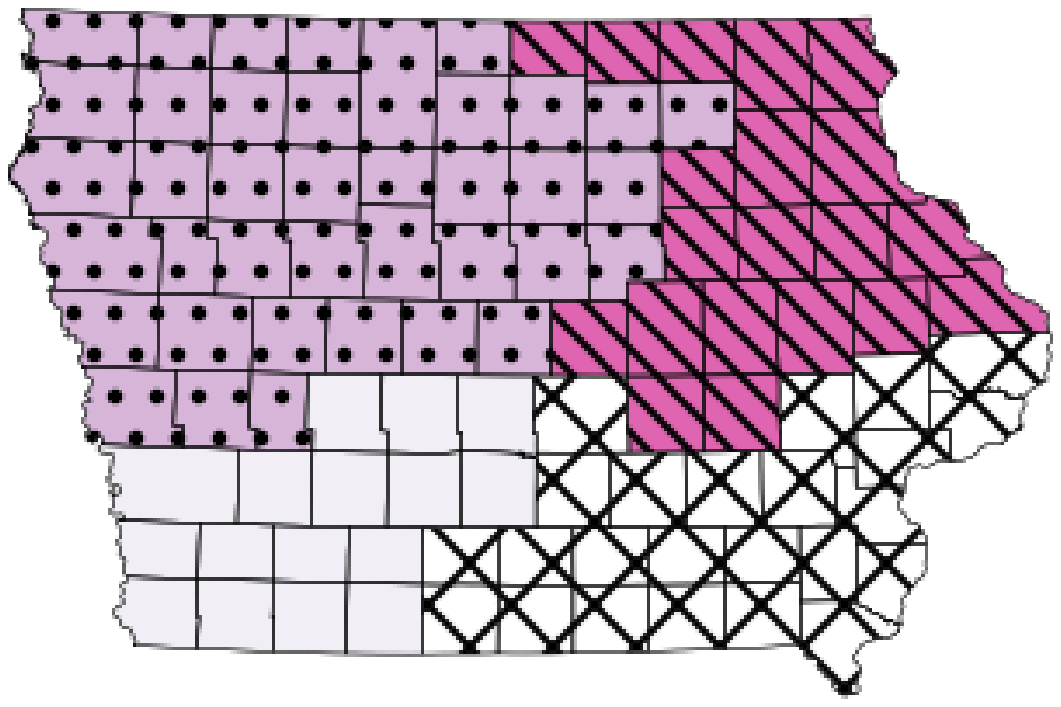}} 
	\subcaptionbox{\label{fig:Swap2}}{\centering \includegraphics[width = 0.48\textwidth]{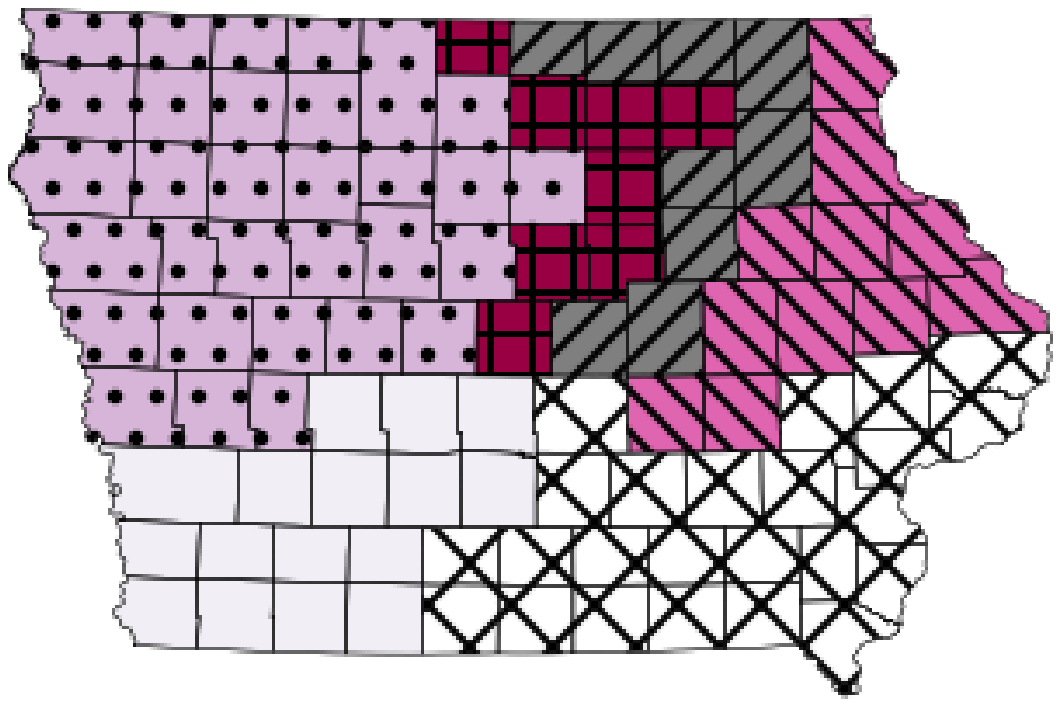}} 
	\subcaptionbox{\label{fig:Swap3}}{\centering \includegraphics[width = 0.48\textwidth]{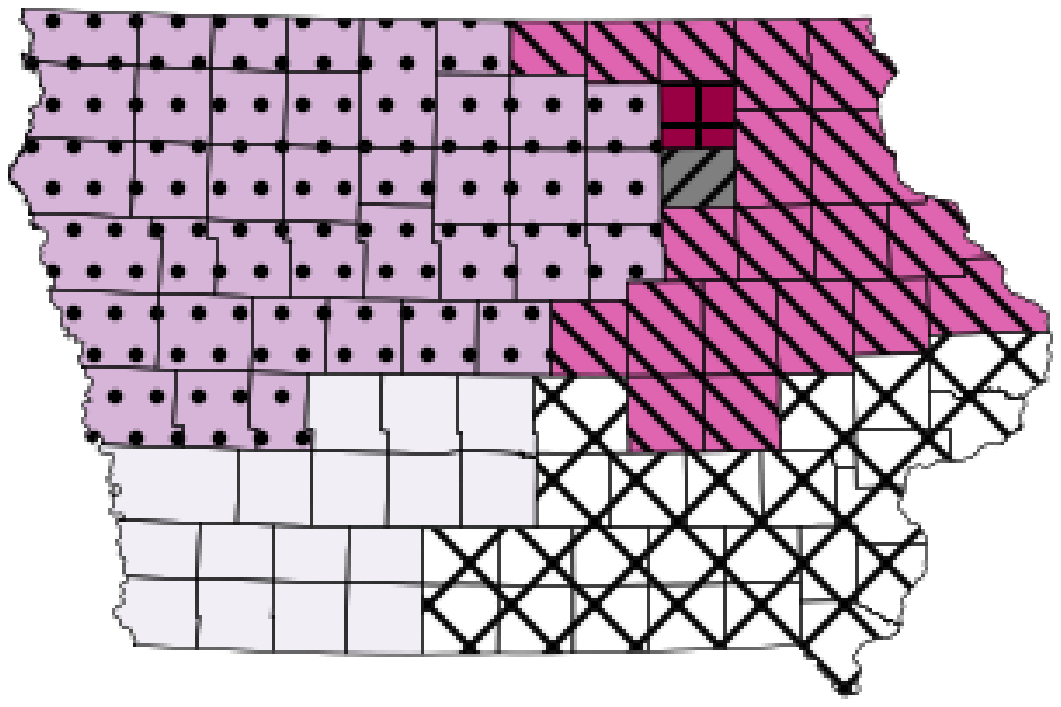}} 
	\subcaptionbox{\label{fig:Swap4}}{\centering \includegraphics[width = 0.48\textwidth]{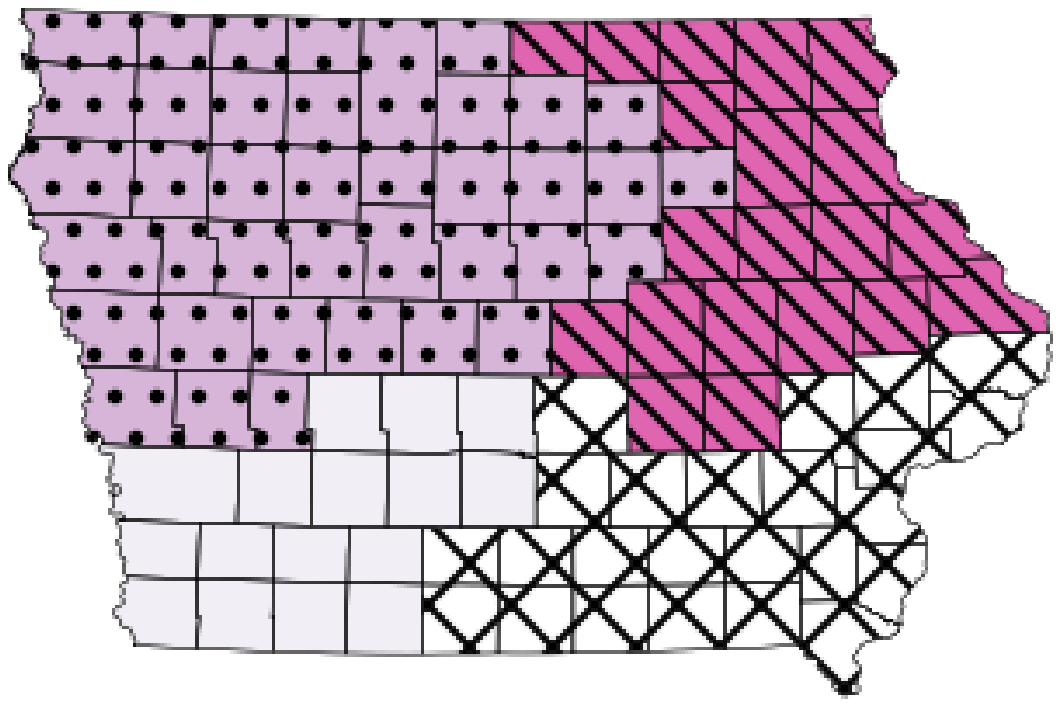}} 
	\caption{Example of Swap neighborhood function. First select Districts that Swap will be performed on, $D_m$: \protect\includegraphics[height =0.25cm]{Fig7.eps} and $D_\ell$: \protect\includegraphics[height = 0.25cm]{Fig8.eps} (\ref{fig:Swap1}), then find pairs of candidate vertices in $D_m$: \protect\includegraphics[height = 0.25cm]{Fig9.eps} and $D_\ell$: \protect\includegraphics[height = 0.25cm]{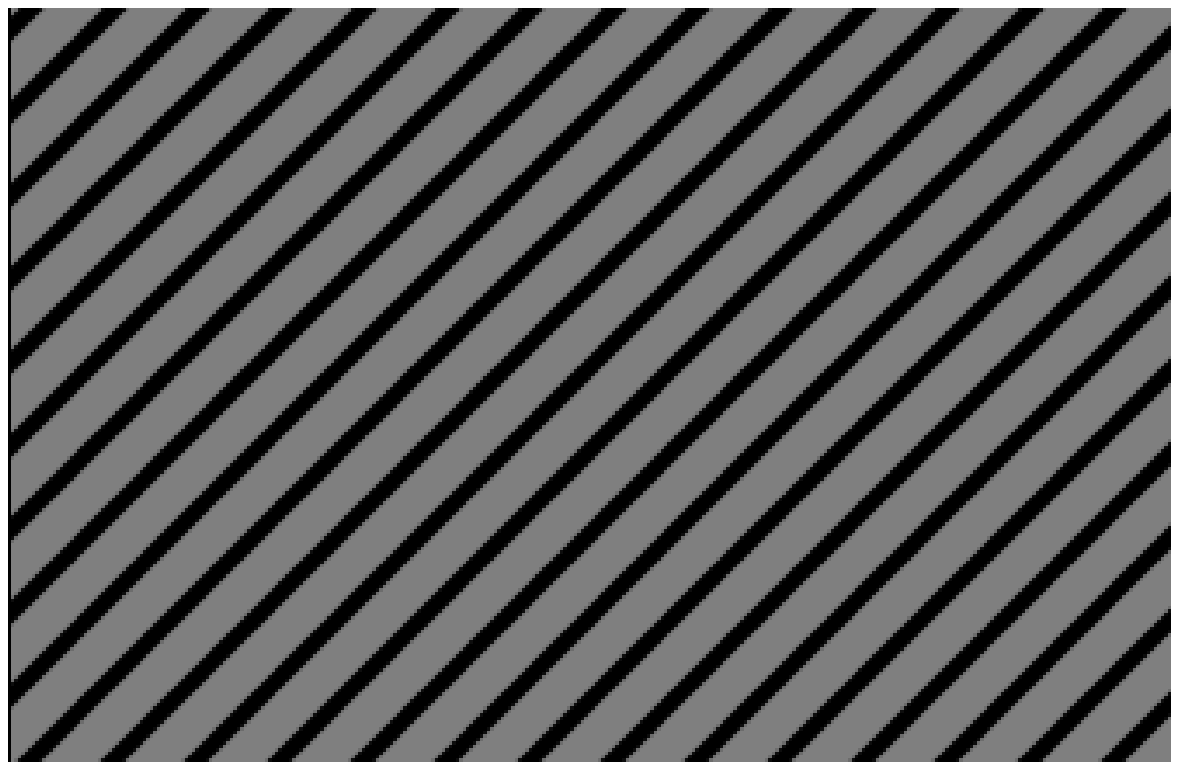} that can be swapped (\ref{fig:Swap2}). Next find candidate vertex pair that minimizes Equation (\ref{LocalSearchMinimizationEquation})  (\ref{fig:Swap3}) and Swap minimizing vertex pair from (\ref{fig:Swap4})
	}
	\label{fig:SwapExample}
\end{figure}

    To construct the districting plan we move $v_i$ to $D_\ell$ and $v_j$ to $D_m$, if $j$ is not the dummy node $-1$. 
    This can be seen in Figure~\ref{fig:FlipExample} for Flip and Figure~\ref{fig:SwapExample} for Swap.

    Starting from this new districting plan, we repeat the process of finding the pair of neighboring districts with the largest disparity, generating neighboring districting plans, finding ($v_i^*$, $v_j^*$) according to Equation (\ref{LocalSearchMinimizationEquation}), and changing districting plans until a districting plan with $\Dev < 1$ is found or until $\LI$ iterations of local search are reached, where $\LI$ is a user defined maximum number of allowed local search iterations.
    
    We noticed that these neighborhood functions could cycle one TB between a few districts, getting stuck in a loop. To get around this we also included a tabu list in the search process, with tabu tenure set as $0.1 \times \LI$. This means that after each successful neighborhood move is made, we store the move in the tabu list, and cannot repeat it for $0.1\times \LI$ iterations. The entire Local Search process is presented in Algorithm \ref{LSalg}.
    
    \begin{algorithm}
    \caption{Local Search}
    \begin{algorithmic}[1]
    \State Choose $\LI$ and $\NF$
    \While{ $\Dev> 1$ and Iterations $<\LI$  }
        \State Store pairs of neighboring districts and population disparity between them
        \While{a valid move has not been found}
        \State Set $\{D_m, D_\ell \}$ as pair of districts with largest population disparity
        \State Find the set of valid, non-tabu moves in $\NH$ according to $\NF$
        \If{there is a valid move}
        \State Make move that minimizes Equation (\ref{LocalSearchMinimizationEquation})
        \State Make move tabu for $0.1\times \LI$ iterations
        \Else
        \State Remove  $\{D_m, D_\ell \}$ from the pairs of neighboring districts
        \EndIf
        \EndWhile
    \State Increment Iteration and compute $\Dev$
    \EndWhile
    \end{algorithmic}
    \label{LSalg}
    \end{algorithm}
    
    \subsection{Coarsening and Uncoarsening}
    \label{section:AlgsCU}
    
    There are $27$ states with counties whose population is larger than $\Pop^*$ for that state. Thus, these states cannot be districted at the county level and instead must be broken down further to the VTD level. This causes the number of TBs to grow very large, and can greatly slow down our method. For example, with Florida there are $67$ counties and $9435$ VTDs. The original graph has $N$ vertices to represent the TBs. The user can provide a parameter $uc_0 \in \left(0, 1\right]$ in order to coarsen $G$ to a graph $G_0$ which has $\lfloor \, uc_0 N \, \rfloor$ vertices. To reduce the number of vertices, the coarsening process randomly selects two neighboring vertices. If the sum of the populations of these two vertices is less than $\Pop^*$, the two neighboring (parent) vertices are combined, resulting in a new (child) vertex.
    
\begin{figure}
	\subcaptionbox{\label{fig:CoarseEx1}}{\centering \includegraphics[width = 0.3\textwidth]{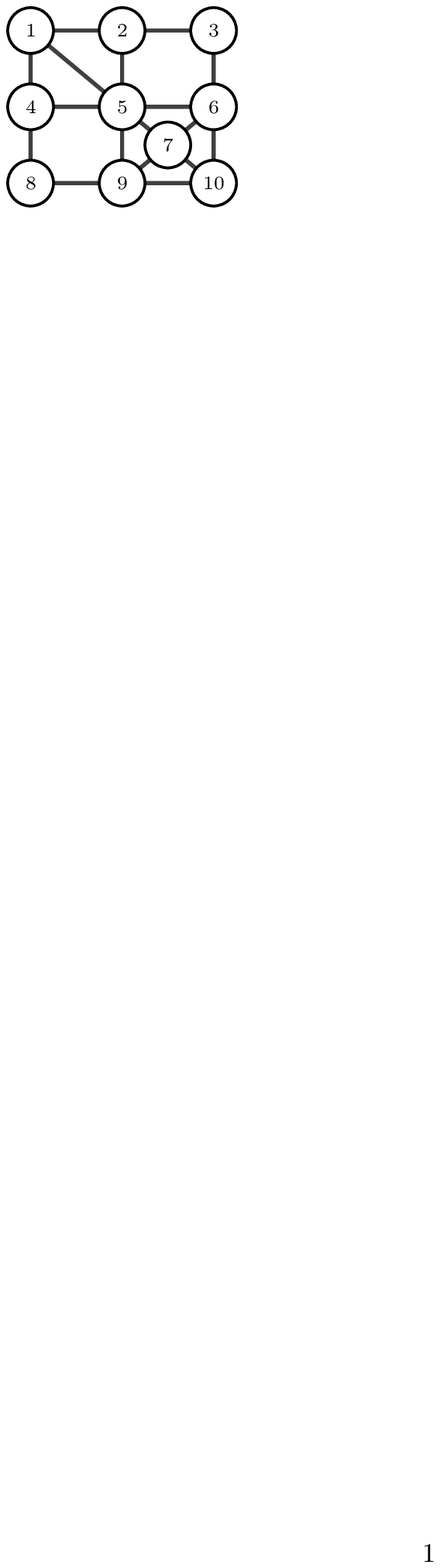}} 
	\subcaptionbox{\label{fig:CoarseEx2}}{\centering \includegraphics[width = 0.3\textwidth]{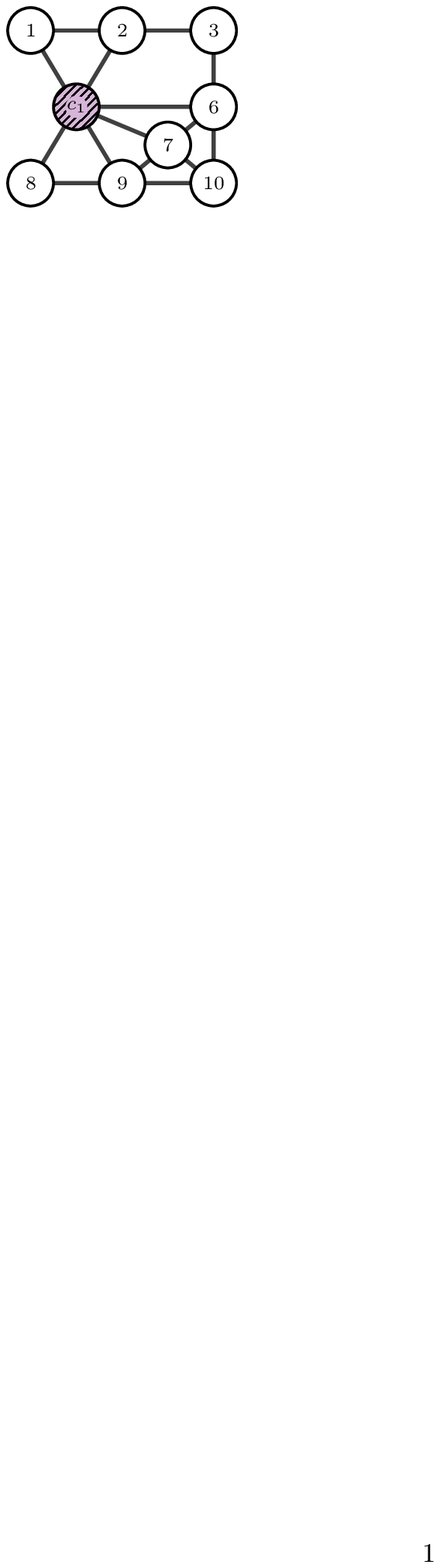}} 
	\subcaptionbox{\label{fig:CoarseEx3}}{\centering \includegraphics[width = 0.3\textwidth]{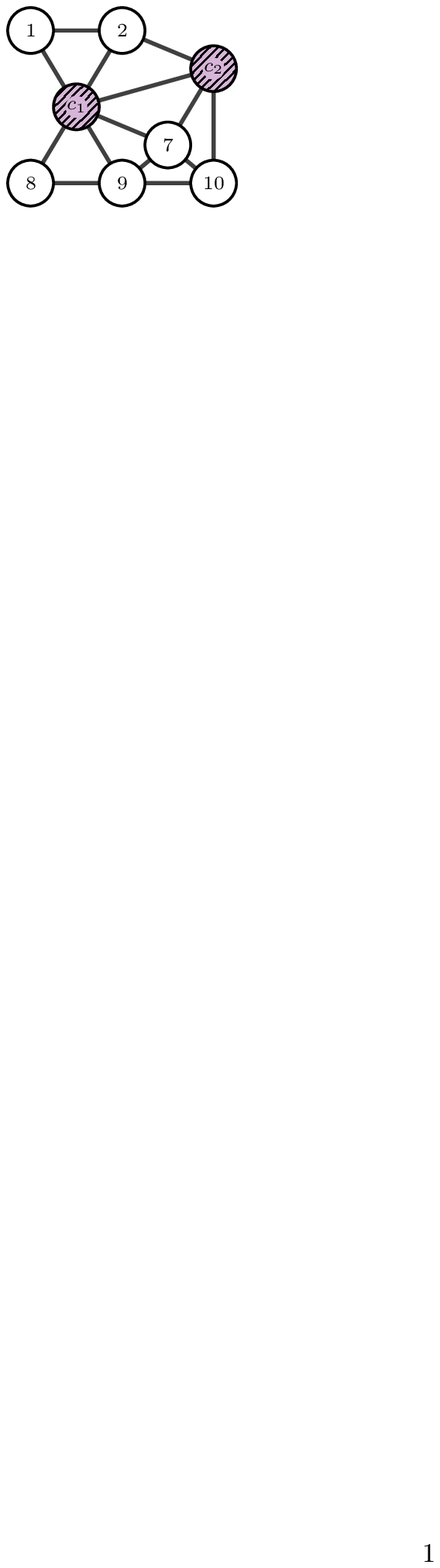}} 
	\subcaptionbox{\label{fig:CoarseEx4}}{\centering \includegraphics[width = 0.3\textwidth]{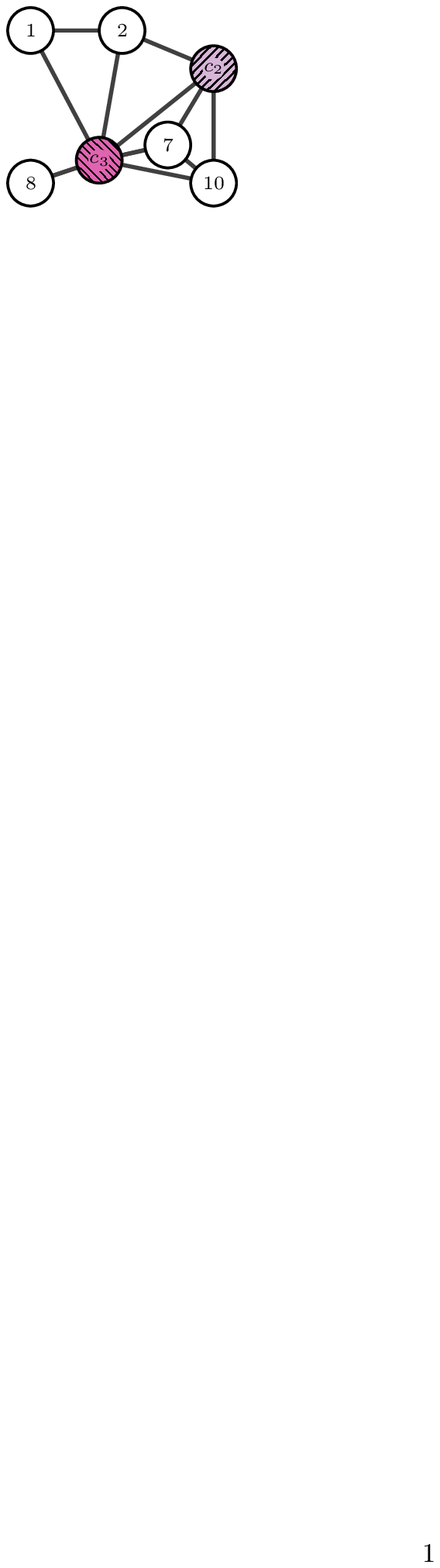}} 
	\subcaptionbox{\label{fig:CoarseEx5}}{\centering \includegraphics[width = 0.3\textwidth]{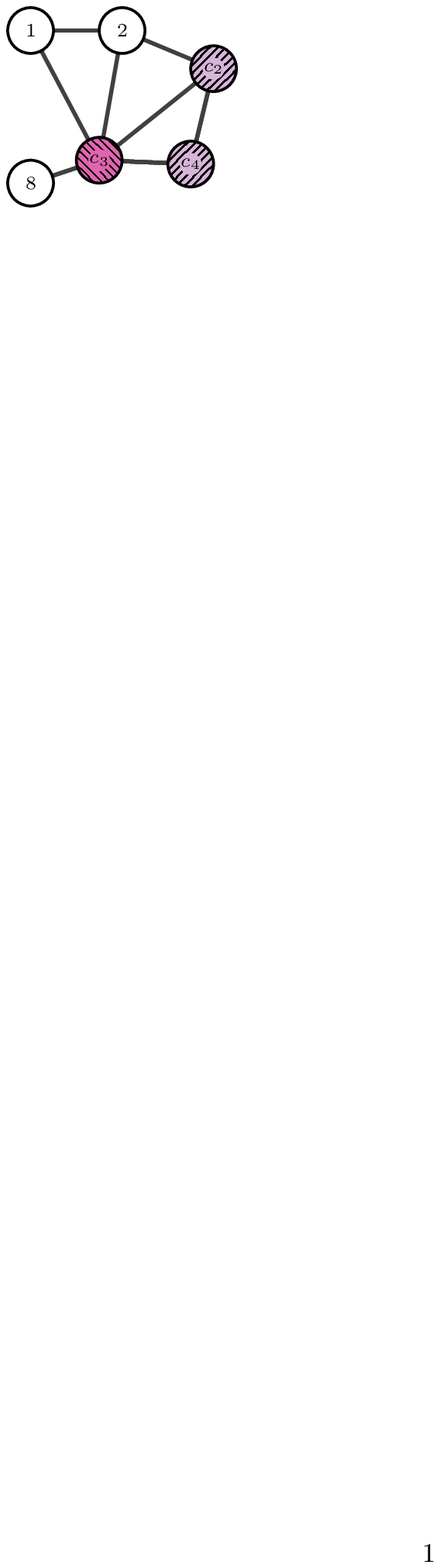}} 
	\caption{Coarsening Example with $uc_{0} = 0.6$. Beginning with graph $G$ which has $N = 10$ vertices (\ref{fig:CoarseEx1}) and $uc_0 = 0.6$, coarsen to $G_0$ with $\lfloor \, uc_0 N \, \rfloor = 6$ vertices. Randomly select pairs of neighboring vertices and combine them until there are only $6$ vertices remaining. First vertices $v_4$ and $v_5$ are combined to form child vertex $c_1$ (\ref{fig:CoarseEx2}), then vertices $v_4$ and $v_5$ are combined to form child vertex $c_1$ (\ref{fig:CoarseEx3}), then vertices $c_1$ and $v_9$ are combined to form child vertex $c_3$ (\ref{fig:CoarseEx4}), lastly $v_7$ and $v_{10}$ are combined to form child vertex $c_4$ and $G_0$ is reached in (\ref{fig:CoarseEx5}). Colors and hatching are used to show coarsened vertices, with darker colors representing more parent vertices, northeast lines representing vertices with two parent vertices, and northwest lines representing vertices with three parent vertices}
	\label{fig:CoarseningExample}
\end{figure}

    The child vertex will inherit edges from the union of the edges of the parent vertices so that the child vertex has the same neighbors as the parent vertices. The population of the child vertex is the sum of the parents' populations. 
    As we coarsen the graph, we store the parent vertices, their populations, edges (or neighbors), and the order they are combined in, to ensure we can reconstruct the original graph, $G$. This coarsening process is repeated until there are only $\lfloor \, uc_0 N \, \rfloor$ vertices in the graph. This process can be seen in Figure~\ref{fig:CoarseningExample}.
  
    In order to rebuild the original graph $G$ from a coarsened graph $G_0$ we use the following uncoarsening process. First we find the most recently created child vertex, remove any edges in $G_0$ that connect to the child vertex and then remove the child vertex from $G_0$. We split the child vertex into it's parent vertices and add both of parent vertices, along with adding any edges that the parent vertices were a part of, into $G_0$. We ensure that the parent vertices retain their original populations and both inherit the district assignment of the child vertex. We repeat this until the graph has $N$ nodes again, this will be the original graph $G$.

\begin{figure}
	\subcaptionbox{\label{fig:UC1}}{\centering \includegraphics[width = 0.3\textwidth]{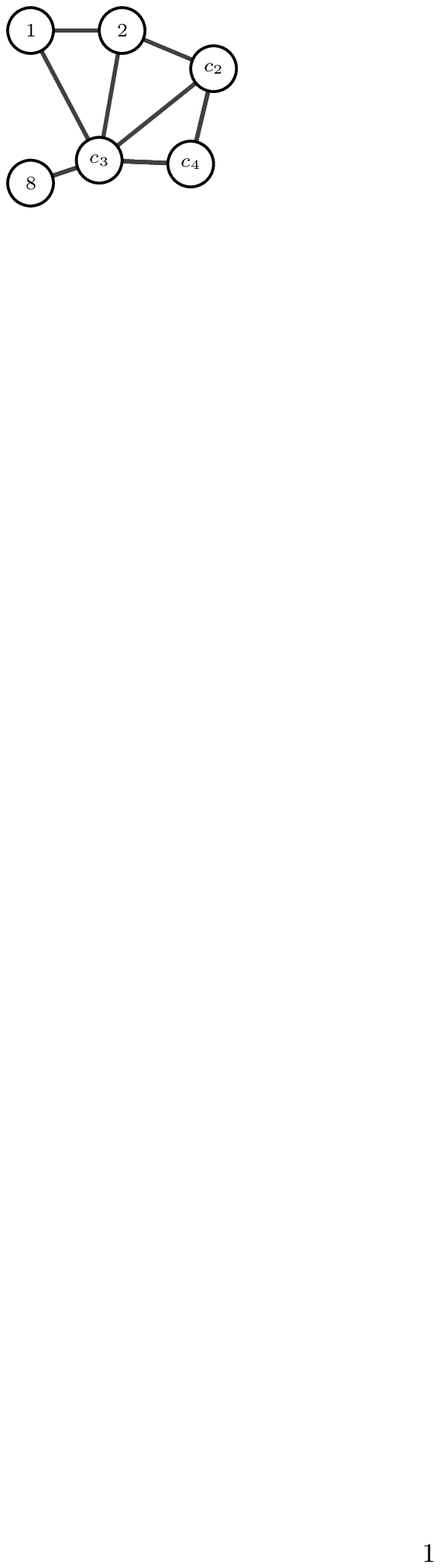}} 
	\subcaptionbox{\label{fig:UC2}}{\centering \includegraphics[width = 0.3\textwidth]{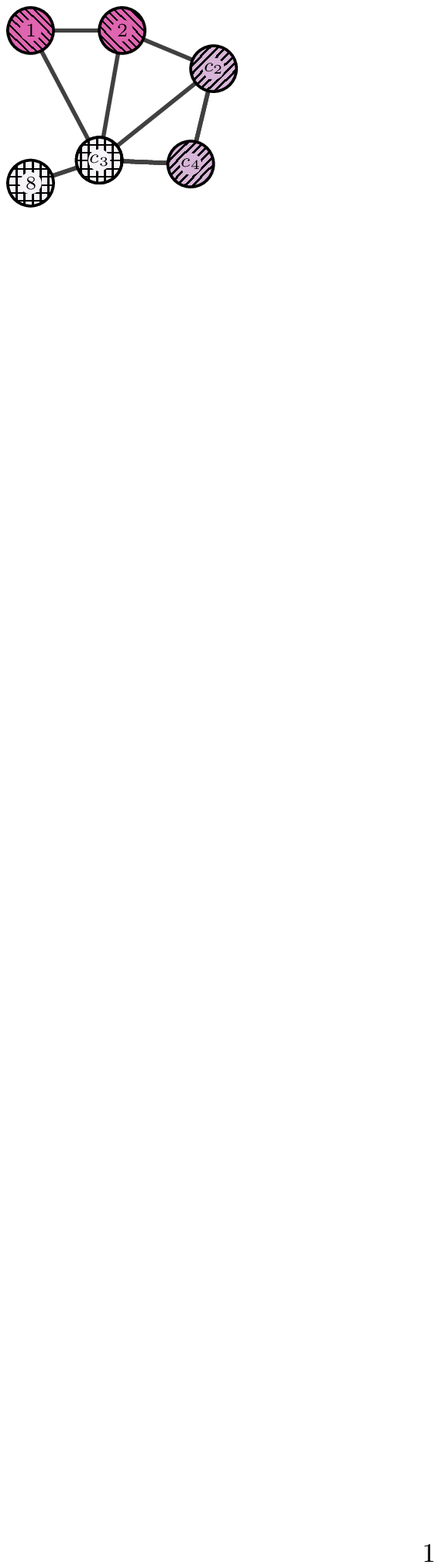}} 
	\subcaptionbox{\label{fig:UC3}}{\centering \includegraphics[width = 0.3\textwidth]{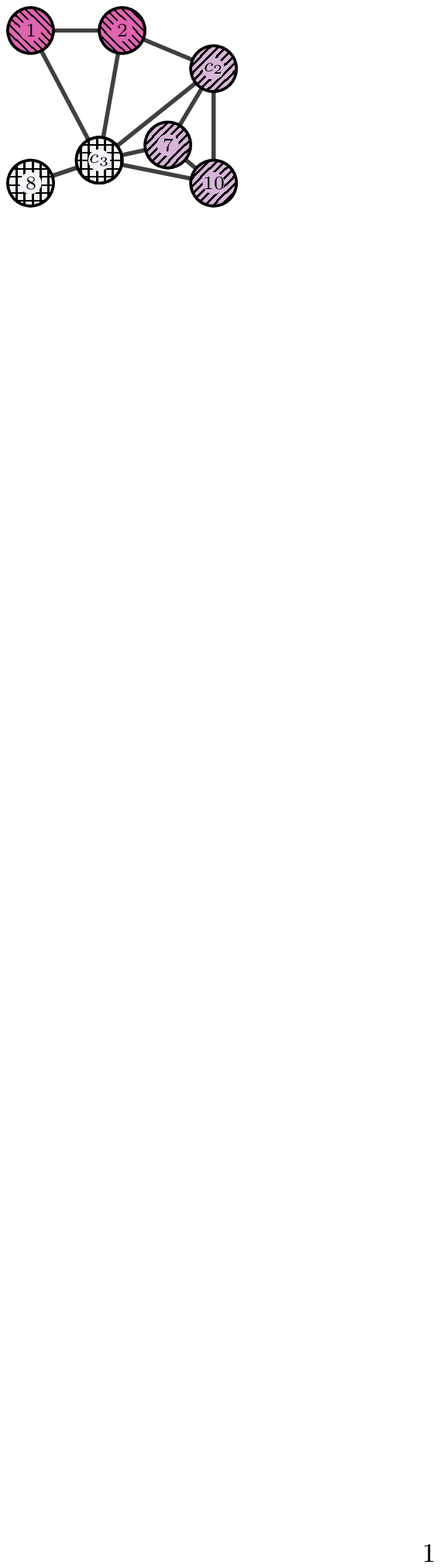}} 
	\subcaptionbox{\label{fig:UC4}}{\centering \includegraphics[width = 0.3\textwidth]{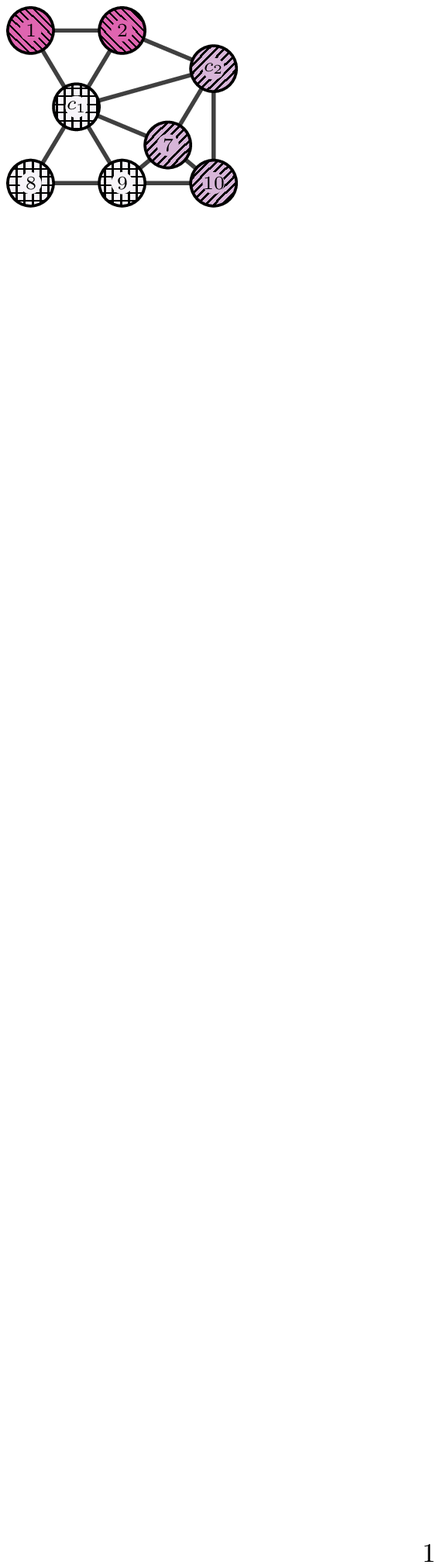}} 
	\subcaptionbox{\label{fig:UC5}}{\centering \includegraphics[width = 0.3\textwidth]{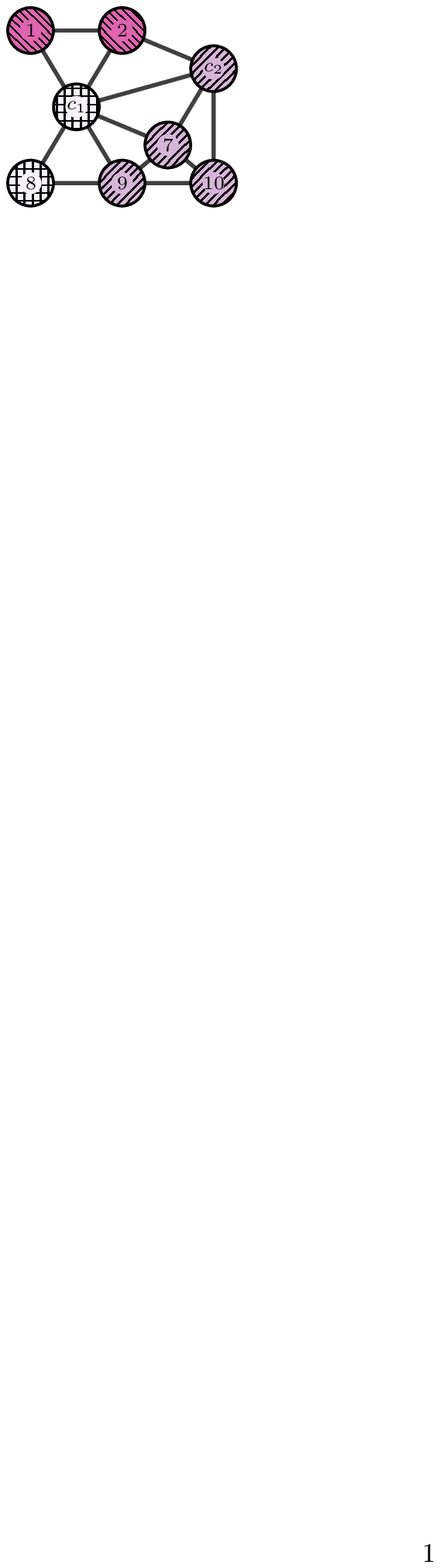}} 
	\subcaptionbox{\label{fig:UC6}}{\centering \includegraphics[width = 0.3\textwidth]{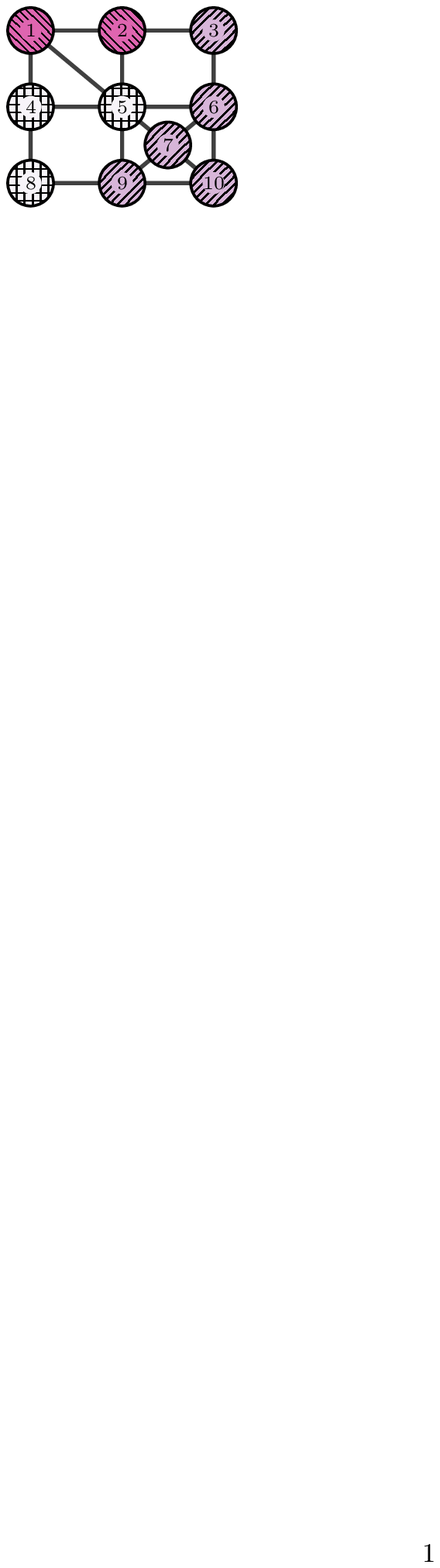}} 
	\caption{Uncoarsening Example with $UC: 0.6, 0.8, 1$ on $G_0$ from Figure~\ref{fig:CoarseningExample} (\ref{fig:UC1}) with $q=2$ uncoarsening steps. From $G_0$ assign each vertex to a district (\ref{fig:UC2}). First uncoarsening step removes most recently created child vertex, $c_4$, from graph and add the parent vertices $v_7$ and $v_{10}$ along with their respective edges (\ref{fig:UC3}). Repeat until there are $\lfloor \, uc_1 N \, \rfloor = \lfloor \, 0.8 \times 10 \, \rfloor = 8 $ vertices and $G_1$ has been found (\ref{fig:UC4}), then reassign vertices to districts through local search (\ref{fig:UC5}).  Repeat until $G_{2}$ has been found and the graph is uncoarsened
	}
	\label{fig:UncoarseningExample}
\end{figure}

    While the coarsening process can help speed up the $k$-medoids process by working with fewer vertices, it has the trade-off of working with larger blocks and does not allow very precise fine tuning, especially in the local search process. To allow for more fine tuning, instead of a single value $uc_0$, the user can provide an uncoarsening schedule with $q$ uncoarsening steps
    $$UC: 0 < uc_0 < uc_1 < \ldots < uc_{q-1} < uc_q = 1.$$
    This will first coarsen $G$ to $G_0$ with $\lfloor \, uc_0 N \, \rfloor$ vertices. Then, instead of uncoarsening to $G$, uncoarsen to $G_1$, a graph with $\lfloor \, uc_1 N \, \rfloor$ vertices. At these partially coarsened steps the local search procedures can be applied to allow for more fine tuning of the solutions. Continue this process of uncoarsening graph $G_i$ to $G_{i+1}$ for $i \in \{0, \ldots, q-1\}$ and doing a local search until $G_q$ is reached. Since $G_q$ is the same as $G$, the original graph will have been reconstructed. One final local search will be run and the uncoarsening process is complete. This process can be seen in Figure~\ref{fig:UncoarseningExample}. Since coarsening may not be needed in every situation, the uncoarsening schedule 
$$UC0: 1$$ 
can be used. This has $uc_0 = uc_q = 1$ and the graph will be coarsened to $\lfloor 1 \times N\rfloor = N$ vertices, which is equivalent to not coarsening.

    \subsection{k-medoids}
    \label{section:AlgsKM}
    
    The final updated $k$-medoids algorithm begins with the selection of an uncoarsening schedule ($UC$),  a maximum number of $k$-medoids iterations ($\MI$), a neighborhood function ($\NF$), and a maximum number of local search iterations ($\LI$). Then $G$ is then coarsened to $G_0$ and $k$ vertices are randomly selected from $G_0$ to be the intial medoids. Each of these new medoids are considered a district, while the remaining vertices are considered unassigned. To assign vertices to districts, the algorithm first determines the district with the smallest population. Then performs one iteration of a Breadth First Search \cite{cormen_leiserson_rivest__1990-BFS} from each of the vertices on the border of this minimum population district to find adjacent vertices. These adjacent vertices are assigned to the district as they are found, provided they are not already assigned to a district and their addition will not cause the district population to exceed $\Pop^*$. If there are no vertices that can be added to the district, the algorithm then moves on to the next smallest district. If no vertices can be added to the adjacent district without exceeding $\Pop^*$ then each unassigned vertex is added to the adjacent district with the smallest population.
    
\begin{figure}
	\subcaptionbox{\label{fig:NewMed1}}{\centering \includegraphics[width = 0.3\textwidth]{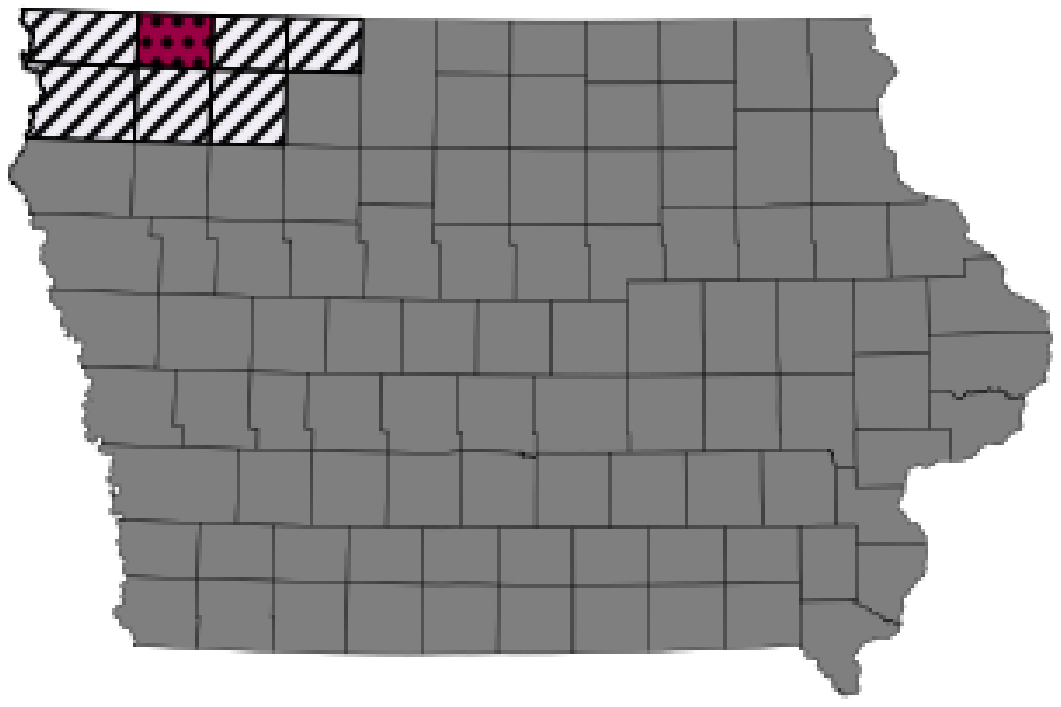}} 
	\subcaptionbox{\label{fig:NewMed2}}{\centering \includegraphics[width = 0.3\textwidth]{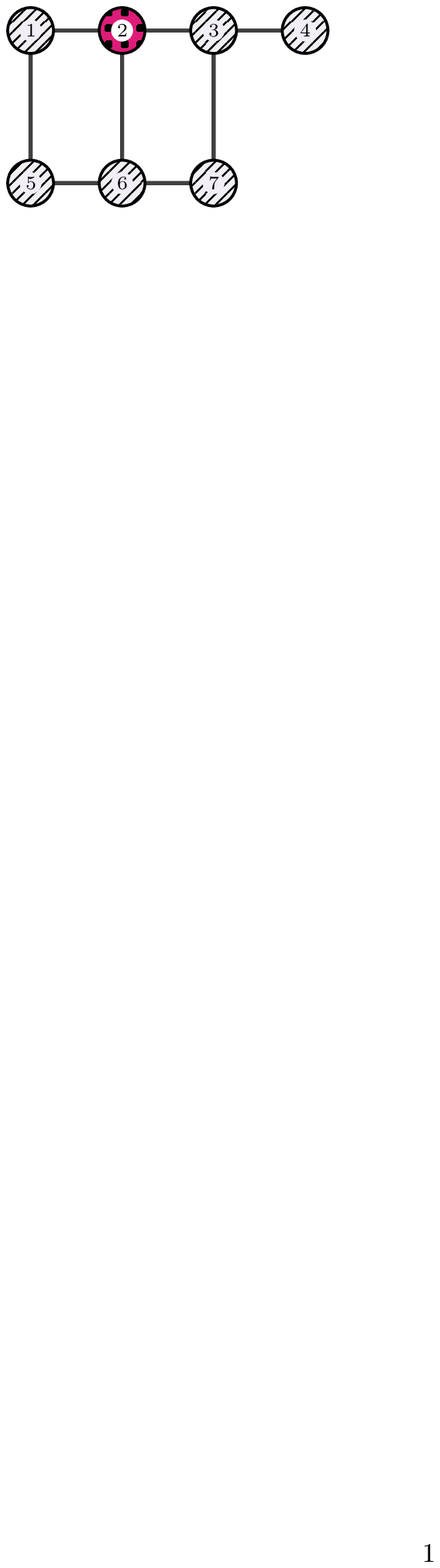}} 
	\subcaptionbox{\label{fig:NewMed3}}{\centering \includegraphics[width = 0.3\textwidth]{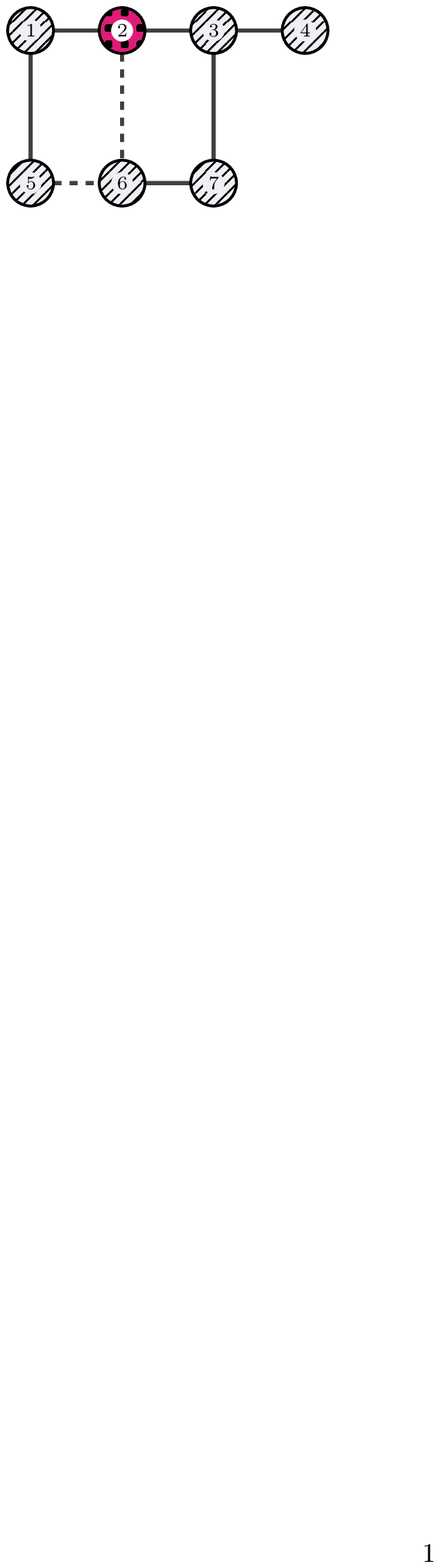}} 
	\subcaptionbox{\label{fig:NewMed4}}{\centering \includegraphics[width = 0.3\textwidth]{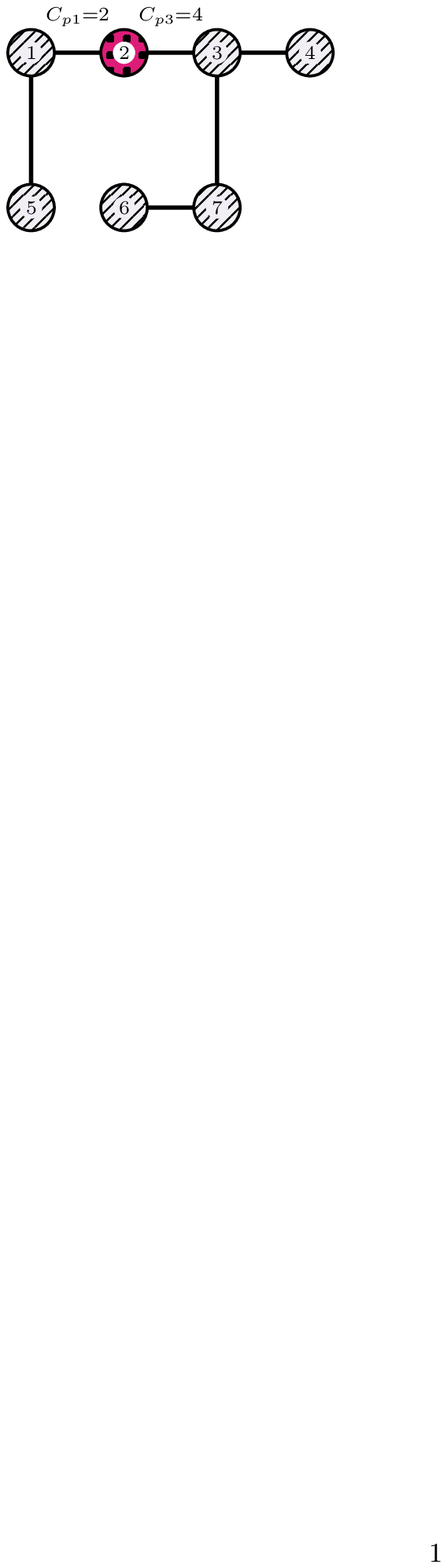}} 
	\subcaptionbox{\label{fig:NewMed5}}{\centering \includegraphics[width = 0.3\textwidth]{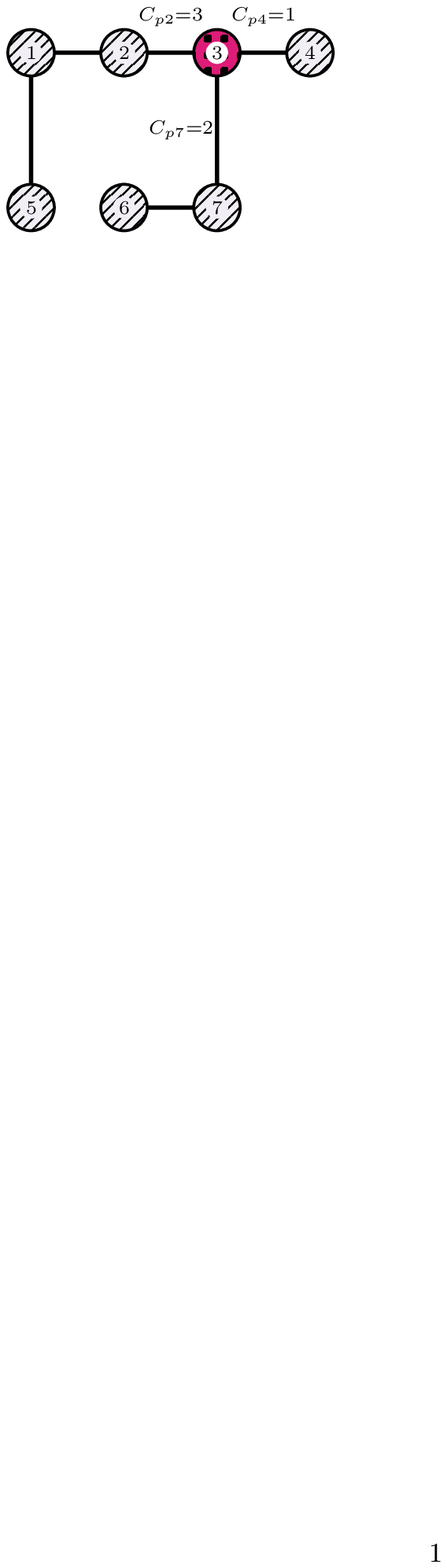}} 
	\subcaptionbox{\label{fig:NewMed6}}{\centering \includegraphics[width = 0.3\textwidth]{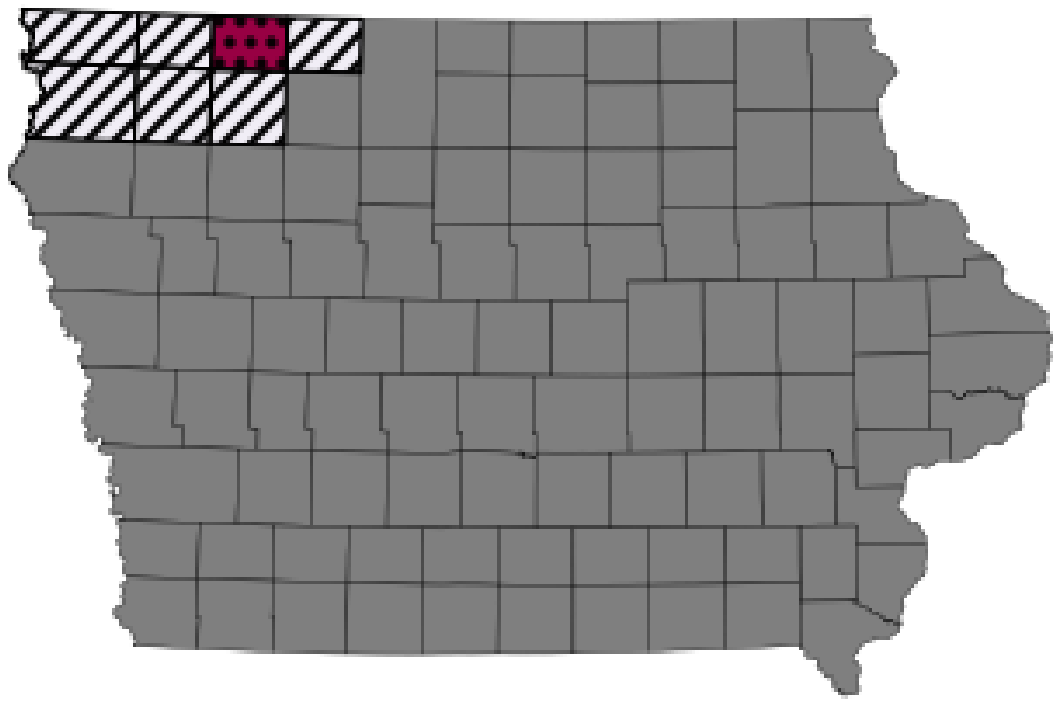}} 
	\caption{Selecting a new medoid example with Iowa. Begin with medoid $v_2$: \protect\includegraphics[height = 0.25cm]{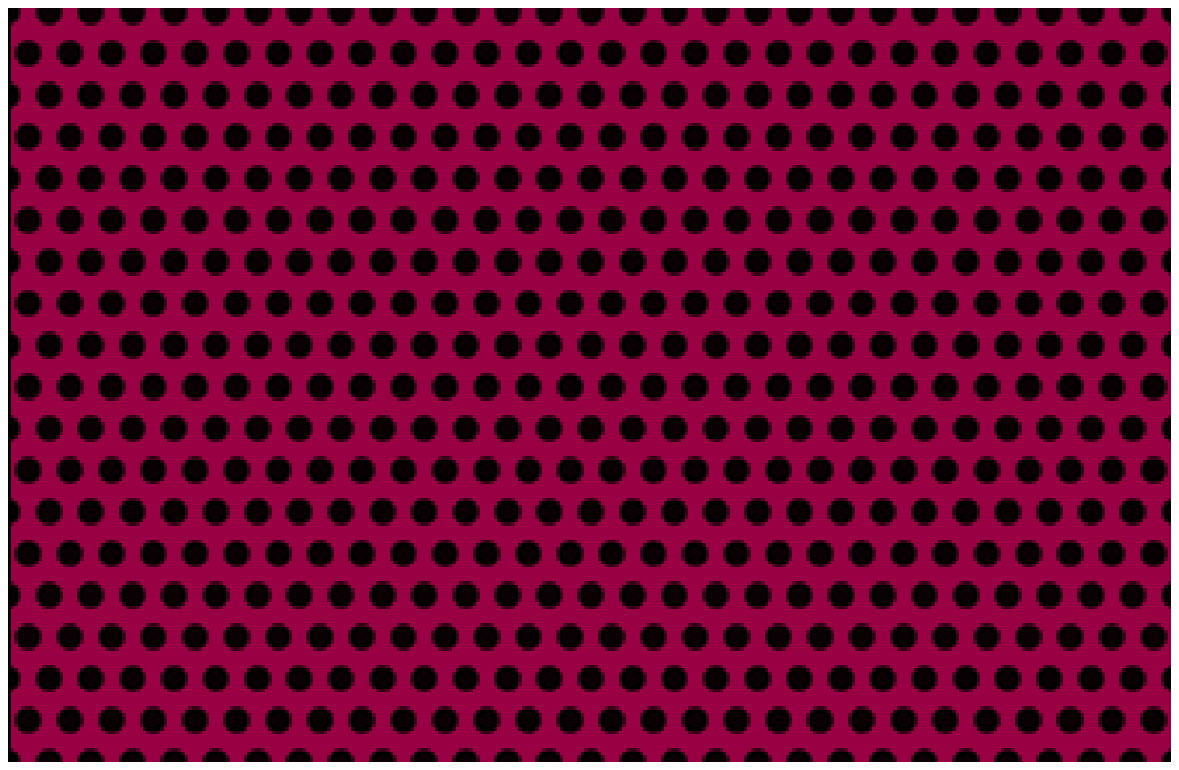} in district: \protect\includegraphics[height = 0.25cm]{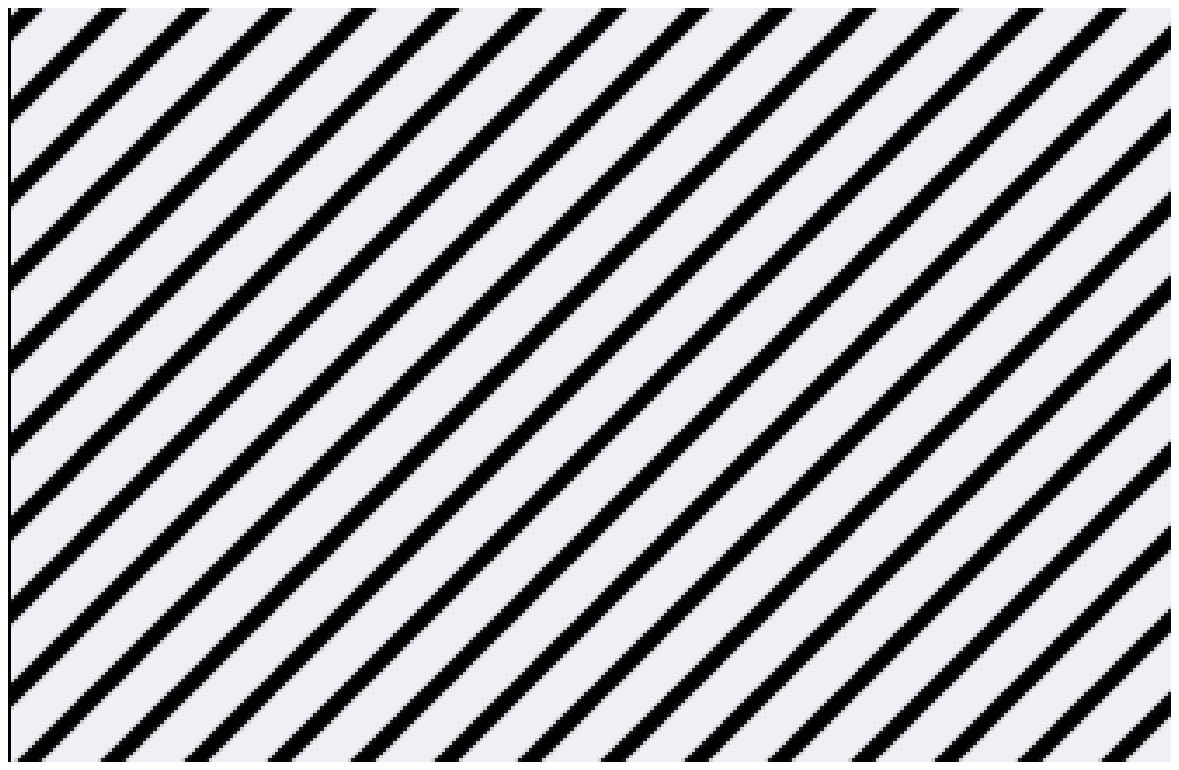} (\ref{fig:NewMed1}). Use graphical representation (\ref{fig:NewMed2}) and apply Broder's Algorithm to randomly select edges to remove to change graph into a tree (\ref{fig:NewMed3}). Dashed lines are used to show the selected edges. Compute length, $C_p$, for each path $p$ beginning at medoid $v_2$ (\ref{fig:NewMed4}). Since $C_{p3} > C_{p1}$ assign first vertex along path $p3$ as new medoid. Compute $C_p$ for each path starting at medoid $v_3$ (\ref{fig:NewMed5}). Since there is no path where $C_{p^*} > \sum_{p\in P_{T,m}\\p^*} C_p$ the medoid does not change and the process terminates. Geographical view of new medoid location (\ref{fig:NewMed6})
	}\label{fig:NewMedoidExample}
\end{figure}

    When all vertices are assigned to districts, the algorithm then determines the medoid, $m$, of each district. This new medoid is selected by first using Broder's Algorithm to randomly cut the edges of a given district until it forms a tree $T$. We let $P_{T,m}$ be the set of paths in $T$ that originate at $m$. For each path $p\in P_{T,m}$ we compute the length, $C_p$, by counting the number of edges along the path. Then we find the path $p^*$ such that $C_{p^*} > \sum_{p\in P_{T,m}\backslash p^*} C_p$, if it exists. 
    Next, we find the vertex neighboring $m$ along $p^*$ and make this vertex the new medoid by reassigning $m$. We repeat until there is no path $p^*$, in that case the current medoid is the new medoid. This new medoid is sensitive to which edges are cut in Broder's Algorithm and since the algorithm determines the new medoid on $T$, the medoid may not appear to be visually centered within the district. An example of the selection process for a new medoids is diagrammed in Figure~\ref{fig:NewMedoidExample}.

    Once a new medoid has been found for each district, compute $\Dev$ for the districting plan. If $\Dev$ is below $1$ or the algorithm has completed $\MI$ iterations, the $k$-medoids process is terminated. Otherwise, all vertices are unassigned from their districts, the iteration variable is incremented, and the district assignment process begins again. Since this process can cause the deviation to increase between iterations, the districting plan with the lowest deviation across all iterations is stored. Once the algorithm has terminated, $\LI$ iterations of local search are performed on the districting plan with the lowest deviation. The user defined uncoarsening schedule is then followed to alternate between uncoarsening and performing a local search until a local search is performed on $G_p$. The complete algorithm is presented in Algorithm~\ref{KMEDSalg}. In addition to the main goal of minimizing the population deviation, this method ensures the generated districting plans are contiguous. In the $k$-medoids process the Breadth First Search construction forces the initial plan to be contiguous as TBs can only be added to adjacent districts. As each new set of medoids is found the Breadth First Search is repeated and the districting plans remain contiguous. In the local search process both Flip and Swap moves maintain contiguity. A Flip will not move a TB if it is an articulation point, thus keeping the districts contiguous. A Swap move will not move a pair $(v_i, v_j)$ if either TBs are articulation points of their respective districts or if they do not share a border with any other TB in the other district, thus keeping the districts contiguous. Since the CMB method is just the best between Flip and Swap moves, it also maintains contiguity. Lastly the uncoarsening process also maintains contiguity, as when a child node is split the parent nodes are assigned to the same district as the child node. The final districting plan criteria is compactness, but this is not something we directly optimize for in our method. Growing the districts from medoids in the $k$-medoids process should keep the districts relatively compact, but we only measure the compactness at the very end of the algorithm after the plan with the lowest deviation has been found. This lesser focus on compactness comes from the Indiana Citizens Redistricting Commission's 2021 public mapping competition where compactness was considered less important than keeping cities and counties intact \cite{ReducedCompactness}.

    \begin{algorithm}
    \caption{$k$-Medoids}
    \begin{algorithmic}[1]
        \State Choose $UC$, $\MI$, $\LI$, and $\NF$
        \State Coarsen the graph to $G_0$
        \State Randomly select $k$ vertices to be the initial medoids
        \While{$\Dev > 1$ and Iterations $< \MI$}
            \While{there are vertices not assigned to a district} 
                \State WorkingDistrict $\leftarrow$ district with minimum population
                \State U $\leftarrow$ Set of all unassigned vertices adjacent to WorkingDistrict
                \State Add vertices of U to WorkingDistrict provided district population does not exceed $\Pop^*$
                \If{no more vertices can be added to any district}
                \State Add each unassigned vertex to adjacent district with minimum population
                \EndIf
            \EndWhile  
            \State Compute medoid of each district, these will be medoids in next iteration
        \EndWhile
        \State Perform $\LI$ iterations of Local Search with $\NF$ on districting plan with smallest deviation
        \For{i = 1, \ldots, q}
        \State Uncoarsen local search output to $G_i$
        \State Perform $\LI$ iterations of local search with $\NF$ on $G_i$
        \EndFor
    \end{algorithmic}
    \label{KMEDSalg}
    \end{algorithm}

\section{Tests and Results}
\label{section:Tests}

	\subsection{Tests}
    For testing we have chosen to use both Iowa and Florida to evaluate our method. Iowa because it is one of the simplest nontrivial cases, having only $99$ TBs and $k=4$ districts. We chose Florida for a more challenging case with $9435$ VTDs and $k=27$ districts. Furthermore, Florida has been used previously in the literature \cite{chen_rodden_2015} and can allow us to more directly compare to the method of Chen and Rodden. To maintain existing state boundaries we added a preprocessing step that combines all VTDs in a county, with total county population less than $\Pop^*$, to a single TB. With this preprocessing step we reduced Florida to $4700$ TBs while keeping $60$ of $67$ counties intact. The population data we used comes from the 2010 census by the United State Census Bureau \cite{PopulationInfoIowa}, \cite{PopulationInfoFlorida} and the shape files we used to make the graphs also come from the United States Census Bureau \cite{ShapeFiles}. As a further note, we will only consider a pair of TBs neighbors if they share a border with nonzero length; TBs that only share a corner point are not considered neighbors. The tests are all run on a single CPU (AMD 2950X at 3.5 GHz) with the algorithms implemented in Python 3.8.5. First, we test parameters to see which allows our model to perform the best. Next, we look at the results with the best parameters. Lastly, we compare to the method used by Chen and Rodden. \cite{chen_rodden_2015}.
    
        \subsubsection{Parameter Selection}
        \label{sec:ParamSel}
        In building the algorithm we evaluated $90$ possible parameter sets where each set consisted of a neighborhood function, a maximum number of local search iterations, and an uncoarsening schedule. We tested three neighborhood functions  ($\NF \in \{ \text{Flip}, \text{Swap}, \text{CMB} \})$, five maximum local search iteration values ($\LI \in \{100, 250, 500, 750, 1000\}$), and six uncoarsening schedules ($UC \in \{UC1, UC2, UC3, UC4, UC5, UC6\}$) defined here:

        \begin{tabular}{cllllllllll} 
            UC1: & 0.3& 0.7&  1 \\
            UC2: & 0.25& 0.5& 0.75& 1\\
            UC3: & 0.3& 0.5& 0.7& 0.9& 1 \\
            UC4: & 0.1& 0.3& 0.5& 0.7& 0.9& 1 \\
            UC5: & 0.3& 0.5& 0.7& 0.9& 0.95& 1 \\
            UC6: &  0.1& 0.2& 0.3& 0.4& 0.5& 0.6& 0.8& 0.9& 1.
        \end{tabular}
        \\
        
Using each of these parameter sets we ran the resulting algorithm $18$ times on Florida, each time with a fixed number of $k$-medoids iterations ($\MI = 100$). The summaries of results for these tests can be seen in Tables~\ref{tab:UC},\ref{tab:NHF},\ref{tab:LSI}. The best found parameters are used in the full tests, if multiple parameters perform well we used both.

	Examining the results for the uncoarsening schedules in Table~\ref{tab:UC}, we can see that $UC6$ outperforms the other uncoarsening schedules with a minimum deviation of $4.28\%$, whereas no other uncoarsening schedules produce districting plans with under $10\%$ deviation. While UC6 is slower than the other uncoarsening schedules, we feel that the lower average deviation and minimum deviation values more than make up for the extra $60$-$120$ seconds of runtime. Even though UC3 is modestly faster than UC5, we claim that UC5 is the second best uncoarsening schedule because it has the next lowest mean deviation after UC6 and approximately matches the minimum deviation of UC3. With this in mind our full tests have been run with uncoarsening schedules UC5 and UC6.

\begin{table}[h]
\begin{minipage}{\textwidth}
\begin{center}
\caption{Uncoarsening Schedule Tests}\label{tab:UC}%
\begin{tabular}{@{}crrrr@{}}
\toprule
 Uncoarsening 	& Mean 		& Min 		& Mean 		& Min \\
 Schedule 	& $\Dev$ 		&  $\Dev$  	& Runtime 	& Runtime \\
 (UC) 		& ($\%$) 		& ($\%$) 		& (s) 		& (s) \\
\midrule
UC1    		& 97.206   	& 13.686 		& 173.903  	& 97.630\\
UC2   		& 92.494   	& 14.207 		& 192.596 	& 107.234 \\
\midrule
UC3    		& 91.716   	& 10.318  		& 200.839  	& 101.113\\
UC4    		& 93.272  		& 10.153 		& 258.924 	& 144.469 \\
\midrule
UC5   		& 86.719  		& 10.897  		& 225.245 	& 107.370 \\
UC6    		& 85.920   	& 4.285  		& 298.919 	& 153.151 \\
\botrule
\end{tabular}
\end{center}
\footnotetext{Uncoarsening Schedule: Uncoarsening Schedule being tested; Mean $\Dev$: mean deviation across all tests; Min $\Dev$: minimum deviation across all tests; Mean Runtime: mean time to coarsen to $G_0$ and run $k$-medoids algorithm - including all uncoarsening and local searches; Min Runtime: minimum time to coarsen to $G_0$ and run $k$-medoids algorithm - including all uncoarsening and local searches}
\end{minipage}
\end{table}
        
        Examining the neighborhood function results in Table~\ref{tab:NHF}, we can see that there is very little difference in the runtimes between the three neighborhood functions. However, there is a clear difference in terms of deviation, with CMB attaining lower deviation in the mean and minimum cases. While we expected CMB to perform better, as it picks between the better of Flip and Swap, the difference between Flip and Swap was surprising. With this in mind our full tests have been run with CMB as the neighborhood function.
        
\begin{table}[h]
\begin{minipage}{\textwidth}
\begin{center}
\caption{Neighborhood Function Tests}\label{tab:NHF}%
\begin{tabular}{@{}crrrr@{}}
\toprule
 Neighborhood	& Mean 		& Min 		& Mean 		& Min \\
 Function 		& $\Dev$ 		&  $\Dev$  	& Runtime 	& Runtime \\
 ($\NF$) 		& ($\%$) 		& ($\%$) 		& (s) 		& (s) \\
\midrule
Flip    		& 89.375  		& 34.173 		& 213.650 	& 97.630\\
\midrule
Swap   		& 135.545   	& 49.255 		& 226.173 	& 98.584 \\
\midrule
CMB    		& 48.744 		& 4.285  		& 235.389  	& 98.182\\
\botrule
\end{tabular}
\end{center}
\footnotetext{Neighborhood Function: Neighborhood Function being tested; Mean $\Dev$: mean deviation across all tests; Min $\Dev$: minimum deviation across all tests; Mean Runtime: mean time to coarsen to $G_0$ and run $k$-medoids algorithm - including all uncoarsening and local searches; Min Runtime: minimum time to coarsen to $G_0$ and run $k$-medoids algorithm - including all uncoarsening and local searches}
\end{minipage}
\end{table}

	Examining the results for maximum number of local search iterations in Table~\ref{tab:LSI}, we see that the $750$ and $1000$ iterations cases outperform the other options. The $750$ and $1000$ iteration cases attain the lowest deviation values of those tested and are only about $30-40$ seconds slower than the fastest cases. While the $500$ iterations case achieves similar deviation values to the $750$ and $1000$ iteration cases, the runtimes are much slower. With this in mind our full tests have been run with $750$ and $1000$ maximum local search iterations.

\begin{table}[h]
\begin{minipage}{\textwidth}
\begin{center}
\caption{Local Search Iterations Tests}\label{tab:LSI}%
\begin{tabular}{@{}crrrr@{}}
\toprule
 Local  Search	& Mean 		& Min 		& Mean 		& Min \\
 Iterations 		& $\Dev$ 		&  $\Dev$  	& Runtime 	& Runtime \\
 ($\LI$) 		& ($\%$) 		& ($\%$) 		& (s) 		& (s) \\
\midrule
100    		& 101.981 	& 19.446 		& 173.893 	& 97.630 \\
\midrule
250   		& 91.352  		& 15.368 		& 229.778 	& 116.340 \\
\midrule
500    		& 89.099 		& 9.995  		& 308.498 	 & 163.629 \\
\midrule
750  			& 86.525   	& 4.285 		& 189.575 	& 125.957 \\
\midrule
1000    		& 85.396 		& 7.446  		& 216.219  	& 138.463\\
\botrule
\end{tabular}
\end{center}
\footnotetext{Local Search Iterations: Maximal number of local search iterations being tested; Mean $\Dev$: mean deviation across all tests; Min $\Dev$: minimum deviation across all tests; Mean Runtime: mean time to coarsen to $G_0$ and run $k$-medoids algorithm - including all uncoarsening and local searches; Min Runtime: minimum time to coarsen to $G_0$ and run $k$-medoids algorithm - including all uncoarsening and local searches}
\end{minipage}
\end{table}

        \subsubsection{Results}
        \label{section:Results}
        We tested the algorithm on Iowa and Florida with slightly different parameter sets. For both states we used $100$ iterations of $k$-medoids ($\MI = 100$), CMB as the neighborhood function ($\NF = CMB$), and two different maximum local search iteration values ($\LI \in \{750, 1000\}$). Since Iowa has so few TBs, the algorithm performs well without any uncoarsening schedule so we chose to not coarsen Iowa at all ($UC = UC0$). Since Florida has many more TBs than Iowa, it produces a much more complicated graph. To handle the complexity of the resultant graph, we used two different uncoarsening schedules ($UC \in \{UC5, UC6 \}$). For each parameter set for each state we ran the resulting algorithm $100$ times - the summary of results for these tests are available in Table~\ref{table:StateResults}. We split the results into two parts, the first being the results from performing the $k$-medoids part of the algorithm on the fully coarsened graph with no local search, the second being the results obtained after performing all of the uncoarsening and local searches.        

\begin{sidewaystable}
\sidewaystablefn
\begin{center}
\begin{minipage}{\textheight}
\caption{Results for $k$-medoids algorithm on Iowa and Florida}\label{table:StateResults}
\begin{tabular*}{\textwidth}{@{\extracolsep{\fill}}crcrrrrrrr@{\extracolsep{\fill}}}
\toprule%
& & & \multicolumn{2}{@{}c@{}}{$k$-medoids} & \multicolumn{3}{@{}c@{}}{Local Search} & Algorithm & Total \\
State & $\LI$ & UC &Mean & Mean & Mean & Min & Mean & Runtime & Runtime \\ \cmidrule{4 - 5} \cmidrule{6 - 8} 
& & & $\Dev$($\%$) & Comp & $\Dev$($\%$)  & $\Dev$($\%$)  & Comp & (s) & (s)\\
\midrule
IA & 750 & UC0 & 1.126 & 0.773 & 0.693 & 0.181 & 0.763 & 0.718 & 3.478 \\
IA & 1000 & UC0 & 1.171 & 0.769  & 0.699 & 0.271 & 0.760 &  0.690 & 3.450 \\
\midrule
FL & 750 & UC5 & 152.221 & 0.707  & 44.059 & 13.915 & 0.626 & 276.853 & 469.880 \\
FL & 750 & UC6 & 163.727 & 0.676  & 40.769 & 12.511 & 0.612  & 328.559 & 512.829 \\
\midrule
FL & 1000 & UC5 & 153.110 & 0.707  & 42.043 & 12.061 & 0.616 & 321.880 & 523.083 \\
FL & 1000 & UC6 & 169.677 & 0.675  & 36.707 & 10.876 & 0.607 & 391.102 & 579.940\\
\botrule
\end{tabular*}
\footnotetext{State: state the tests were done on; $\LI$: number of local search iterations; UC: uncoarsening schedule followed; $k$-medoids Mean $\Dev$: mean deviation of the algorithm without any local search; $k$-medoids Mean Comp: mean compactness score of the algorithm without any local search; Local Search Mean $\Dev$: mean deviation after $\LI$ iterations of local search; Local Search Mean Comp: mean compactness after $\LI$ iterations of local search; Local Search Min $\Dev$: minimum deviation across all tests with given initial conditions; Algorithm Runtime: mean time to coarsen to $G_0$ and run $k$-medoids algorithm - including all uncoarsening and local searches; Total Runtime: Algorithm Runtime plus mean time to read in data and summarize results}
\end{minipage}
\end{center}
\end{sidewaystable}

From Table~\ref{table:StateResults} we can see that our method works very well for Iowa. Without local search the $k$-medoids algorithm is able to attain a mean deviation of about $1.1\%$, just over the legally allowed $1\%$ deviation. Including local search decreases the mean deviation to about $0.69\%$, just below the legal limit. The local search also slightly decreases in the mean compactness value (only a difference of  $\sim 0.01$). There is hardly any difference between the two parameter settings for Iowa, with the $750$ iterations case finding slightly lower deviation values, slightly higher compactness scores, and running a bit slower than the $1000$ iterations case.

Our method is less consistent for the more complicated case of Florida. While including local search decreases the mean deviation value by more than $100\%$, the minimum deviation values found were still higher than $1\%$. As with the Iowa case, we can see that adding local search decreases the mean compactness score. Comparing the parameter settings we can see that the final districting plans produced using UC6 tend to have lower deviations than final plans produced using UC5. We can also see that districting plans produced with $1000$ iterations of local search tend to have lower deviations than plans produced with $750$ iterations of local search. Here we also note the random nature of this algorithm, as the minimum deviation in Table~\ref{table:StateResults} is $10.876\%$ for UC6 with $1000$ iterations of local search, while the minimum deviation found in the Section~\ref{sec:ParamSel} parameter setting tests was $4.285\%$ for UC6 with $750$ iterations of local search.

With our method consistently finding districting plans with deviation below $1\%$ on Iowa, we added in an extra step to find more good districting plans over the course of a single run. Anytime the algorithm finds a districting plan that has a deviation less than $5\%$ during the $k$-medoids phase, we save it. Then we perform the local search step on all of these additional districting plans as well. The results of these tests are in Table \ref{table:AdditionalIowa}. With this additional step we are able to find another nine districting plans that have a mean deviation of about $0.6\%$ each time we run our algorithm. Moreover, these extra districting plans come at a cost of about $2.77$ seconds each, slightly lower than the approximately $3.46$ seconds for the first districting plan. Overall this leads to finding about ten districting plans below $1\%$ deviation in just under thirty seconds. 

 \begin{table}
\begin{minipage}{\textwidth}
\caption{Additional District Results} \label{table:AdditionalIowa}
\begin{center}
 \begin{tabular}{crrrrrrr}
 \toprule
  	  & 		& 	 	& Mean 	     	&Mean           	& Mean	& Mean		& Total \\
 State & $\LI$   & UC 	&  Additional 	& $\Dev$        	&  Comp	& Runtime		& Runtime\\
 	  &    		&  	 	&  Districts    	& Per ($\%$)  	&  Per	& Per (s)		& (s) \\
 \hline
 IA 	  & 750 	& UC0 	& 9.290 		& 0.667  		& 0.755 	& 2.775		& 28.577\\ 
 IA 	  & 1000 	& UC0 	& 8.830 		& 0.671		& 0.775	& 2.775		& 28.642\\
 \botrule
 \end{tabular}
 \end{center}
\footnotetext{State: state tests were done on; $\LI$: number of local search iterations; UC: uncoarsening schedule followed; Mean Additional Districts: mean number of districting plans with deviation below $5\%$ found during $k$-medoids process; Mean $\Dev$ Per: mean deviation of each additional districting plan after local search is performed; Mean Comp Per: mean of compactness value for each additional districting plan after local search is performed; Mean Runtime Per: mean time to perform a local search on each additional districting plan and summarize results; Total Runtime: mean time to read in data, run $k$-medoids algorithm, perform local searches on additional districting plans and summarize all results}
\end{minipage}
 \end{table}

    \subsection{Comparison to Alternative Method}
    Next we wanted a direct comparison to a different, but similar, method. We implemented our own version of the algorithm described by Chen and Rodden in \cite{chen_rodden_2015}. We ran tests on this method to directly compare to the best versions of $k$-medoids. Since coarsening allows the local search to occur many times, directly comparing local search iterations between the two methods would be an uneven number of searches. To make up for this difference, we allowed Chen and Roddens' method take $\LI\times \lvert UC \rvert$ iterations. This way they each use the same number of total local search iterations. We compare the methods using the same tests from Section \ref{section:Results}.
   
\begin{sidewaystable}
\sidewaystablefn
\begin{center}
\begin{minipage}{\textheight}
\caption{Comparison of $k$-medoids Algorithm with Chen and Rodden method}\label{table:CRCompResults}
\begin{tabular*}{\textwidth}{@{\extracolsep{\fill}}ccrcrrrrr@{\extracolsep{\fill}}}
\toprule
& & &  &\multicolumn{3}{@{}c@{}}{Local Search}  &Algorithm  & Total \\ 
Method & State & $\LI$ & UC &Mean & Min & Mean & Runtime & Runtime \\ \cmidrule{5-7} 
& & &  & $\Dev$($\%$) & $\Dev$($\%$) &  Comp & (s) & (s)\\
\midrule
KMED & IA & 750 & UC0 & 0.693 & 0.181 & 0.763 & 0.718 & 3.478 \\
CR & IA & 750 & UC0 & 7.162 & 0.426 & 0.782 & 28.155 & 29.663 \\
\midrule
KMED & IA & 1000 & UC0 & 0.699 & 0.271 & 0.760 & 0.690 & 3.450 \\
CR & IA & 1000 & UC0 & 6.967 & 0.104 & 0.784 & 34.914 & 36.388 \\
\midrule
KMED & FL & 750 & UC5 & 44.059 & 13.915 & 0.626 & 276.853 & 469.880 \\
CR & FL & 4500 & UC5 & 143.884 & 65.195 & 0.692 & 706.139 & 728.052 \\
\midrule
KMED & FL & 750 & UC6 & 40.769 & 12.511 & 0.612 & 328.559 & 512.829 \\
CR & FL & 7500 & UC6 & 147.769 & 41.759 & 0.697 & 1080.773 & 1102.027 \\
\midrule
KMED & FL & 1000 & UC5 & 42.043 & 12.061 & 0.616 & 321.880 & 523.083 \\
CR & FL & 6000 & UC5 & 147.096 & 50.079 & 0.694 & 774.962 & 797.021 \\
\midrule
KMED & FL & 1000 & UC6 & 36.707 & 10.876 & 0.607 & 391.102 & 579.940 \\
CR & FL & 10000 & UC6 & 143.654 & 51.167 & 0.699 & 1201.172 & 1223.036\\
\botrule
\end{tabular*}
\footnotetext{Method: method used $k$-medoids algorithm (KMEDS) or the Chen and Rodden method (CR); State: state the tests were done on; $\LI$: number of local search iterations; UC: uncoarsening schedule followed; Local Search Mean $\Dev$: mean deviation after $\LI$ iterations of local search; Local Search Mean Comp: mean compactness after $\LI$ iterations of local search; Local Search Min $\Dev$: minimum deviation across all tests with given initial conditions; Algorithm Runtime: average time to run the algorithm (including uncoarsening and local search time if used); Total Runtime: Algorithm Runtime plus average time to read in data and summarize results}
\end{minipage}
\end{center}
\end{sidewaystable}

    From Table \ref{table:CRCompResults} we can see that our $k$-medoids method produces districts with lower deviation than the other method in the same number of iterations. For example, in the Iowa tests our method has a mean deviation below $1\%$, while the other method has a mean deviation of about $7\%$. In the Florida tests our method has a mean deviation of about $40\%$, while the other method has a mean deviation of about $145\%$. Since our method runs faster than the Chen and Rodden method ($\sim 20-30$ seconds for Iowa and $\sim300-500$ seconds for Florida), our method is able to find more districting plans faster. With our method also tending to have lower deviation values, our method is more likely to find a good districting plan in less time. Furthermore, if we account for the extra districting plans found with by local searching on any plan found with $\Dev < 5\%$, our method is able to find about ten districting plans with $\Dev < 1\%$ in the approximately the same time it takes the other method to finish. However, the Chen and Rodden method does produce districts that are slightly more compact than ours, with the mean compactness values being $0.02$ to $0.09$ higher. This indicates to us that our method outperforms the Chen and Rodden method. 
    
    \subsection{Current Districting Plans}
    
For a baseline comparison we look at what the official districting plans for Iowa and Florida are from the $2010$ Census. 

\begin{figure}
	\subcaptionbox{\label{fig:IowaCurrent}}{\centering \includegraphics[width = 0.45\textwidth]{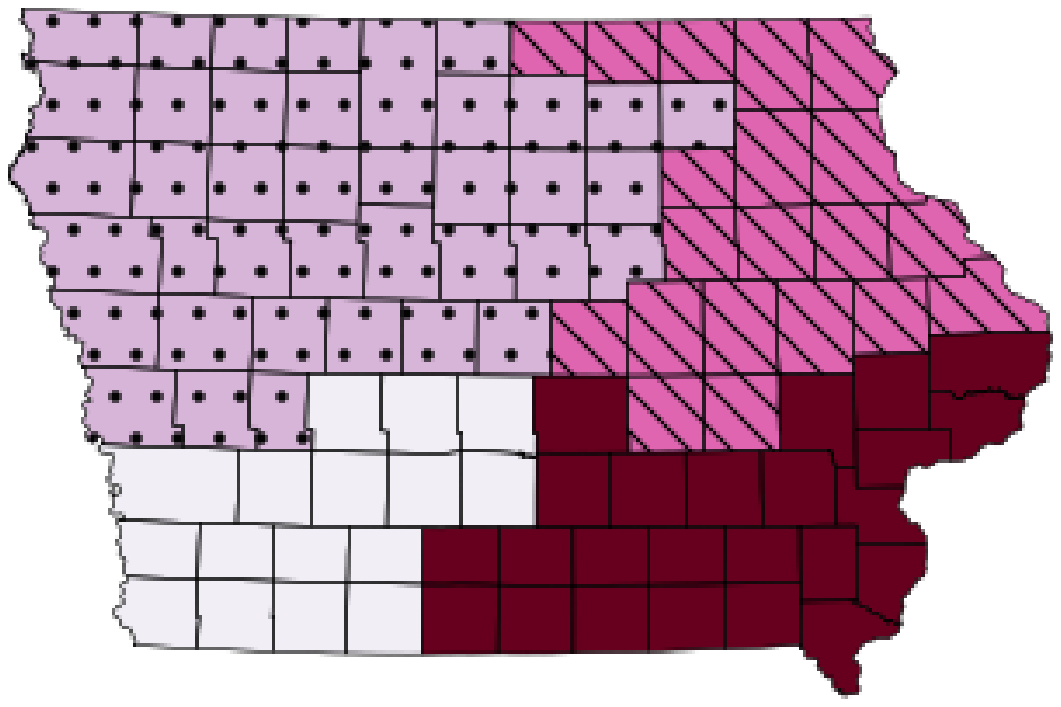}} 
	\subcaptionbox{\label{fig:FloridaCurrent}}{\centering \includegraphics[width = 0.4\textwidth]{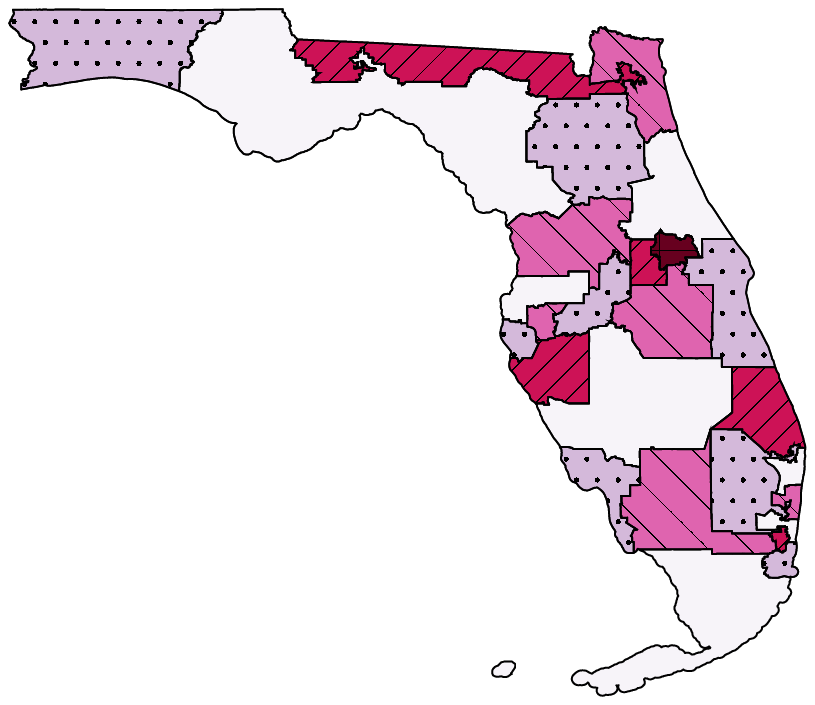}} 
	\caption{Districting plans used by Iowa (\ref{fig:IowaCurrent}) and Florida (\ref{fig:FloridaCurrent}) in 2021
	}
	\label{fig:TrueDP}
\end{figure}

The official plan Iowa used in 2021 can be seen in Figure \ref{fig:IowaCurrent}. It has a maximum population deviation of $0.005\%$ and compactness of $0.78$. Iowa law requires all districting plans to have deviation less than $1\%$ and ours are able to achieve the legal requirement, but with deviation slightly higher than the official Iowa districting plan. Additionally, our districting plans have very comparable compactness values, averaging about $0.76$.

The official plan for Florida used in 2021 can be seen in Figure \ref{fig:FloridaCurrent}. It has a maximum population deviation of $0.0007\%$ and compactness of $0.767$, which equates to districts with populations that differed by at most $1$ person. However, this was done using a finer grain TB than we used, breaking the state into Census Blocks. For Florida there are $488553$ Census Blocks, allowing for much more fine tuning at the cost of less intuitive borders. This plan has a deviation at least $4.2847$ lower than the best plan our method has created and has a compactness value about $0.1$ higher than the plans produced by our method.

\section{Conclusions and Future Work}
Our method works very well on the simple case of Iowa, and performs reasonably on the more challenging case of Florida. In both cases it outperforms the method used by Chen and Rodden in terms of both time and deviation achieved. Our method is able to combine many of the smaller ideas used throughout the years into one method. We also introduce a seemingly novel method to build initial districts by allowing the smallest district to expand next at each step. 

Further work needs to be done working on the remaining states and with the updated census data. Furthermore, we could add additional methods to potentially improve our algorithm. Examples would be the addition of multiple object moves in the local search step from \cite{jin_2017}, the implementation of geographs for faster local search \cite{Article3-GeoGraphs,Article5-GeographPractice}, forcing the local search to only make moves that strictly decrease the deviation, different selection criteria for which pair of districts to use in the local search phase.

We could also add in a specific piece to create minority-majority districts in states where they are required. An approach for this would be to merge neighboring TBs with high minority populations before any clustering is done. This would force the minority communities of interest to stay intact and could lead to the minority-majority districts.

\section*{Declarations}

\subsection*{Competing Interests}
	On behalf of all authors, the corresponding author states that there is no conflict of interest.
	
	
\subsection*{Funding}
	We wish to acknowledge the support of the National Science Foundation for Anthony Pizzimenti through their grant NSF-1407216.

\subsection*{Availability of data and materials}
	The census data analyzed during the current study are available from the census bureau's website. The Florida data is at https://data.census.gov/cedsci/table?text=P1\&g=0400000US12\%24700000\\\&y=2010\&tid=DECENNIALPL2010.P1 and the Iowa data is at https://data.census.gov/cedsci/table?text=P1\&g=0400000US19\%24700000\\\&y=2010\&tid=DECENNIALPL2010.P1. Alternatively, the data can be reached from https://data.census.gov, searching for table P1, selecting the year 2010, selecting Geography tab $\rightarrow$ Voting District $\rightarrow$  (Iowa or Florida) $\rightarrow$  All Voting Districts (VTD).

	The Shape files used during the current study are available from https://www2.census.gov/geo/tiger/TIGER2012/VTD/ .
	
	The code and results data that support the findings of this study are available from the corresponding author upon request.

\bibliography{kmedoids_Grove.bib}

\end{document}